\documentclass[aos,MSNbibl,seceqn,dvips]{arximspdf}
\usepackage[mathcal]{euscript}
\usepackage{graphicx}
\doi{10.1214/13-AOS1188} 
\volume{42}
\issue{1}
\pubyear{2014}
\firstpage{255}
\lastpage{284}

\makeatletter
\newproclaim{definition}{Definition}
\newproclaim{remark}{Remark}
\newtheorem{theorem}{Theorem}

\newtheorem{proposition}{Proposition}
\newtheorem{corollary}{Corollary}
\newcommand{\R}[1]{\mathbb{R}^{#1}}
\newcommand{\eps}{\varepsilon}
\newcommand{\ind}{\mathbb{I}}
\newcommand{\pr}{\mathbb{P}}
\newcommand{\A}{\mathcal{A}}
\newcommand{\F}{\mathfrak{F}}
\newcommand{\M}{\mathcal{M}}
\newcommand{\C}{\mathcal{C}}
\makeatother

\begin{document}
\begin{frontmatter}

\title{On Poincar\'e cone property}
\runtitle{On Poincar\'e cone property}

\begin{aug}
\author[a]{\fnms{Alejandro} \snm{Cholaquidis}\ead[label=e1]{acholaquidis@cmat.edu.uy}},
\author[b]{\fnms{Antonio} \snm{Cuevas}\corref{}\ead[label=e2]{antonio.cuevas@uam.es}\ead[label=u1,url]{http://www.uam.es/antonio.cuevas}\thanksref{t1}}
\and
\author[a]{\fnms{Ricardo}~\snm{Fraiman}\ead[label=e3]{rfraiman@cmat.edu.uy}\thanksref{t1}}
\thankstext{t1}{Supported in part by Spanish Grant MTM2010-17366.}
\runauthor{A. Cholaquidis, A. Cuevas and R. Fraiman}
\affiliation{Universidad de la Rep\'ublica,
Universidad Aut\'onoma de Madrid
and Universidad de la Rep\'ublica}
\address[a]{A. Cholaquidis\\
R. Fraiman\\
Centro de Matem\'atica\\
Universidad de la Rep\'ublica\\
11400-Montevideo\\
Uruguay\\
\printead{e1}\\
\phantom{E-mail:\ }\printead*{e3}}

\address[b]{A. Cuevas\\
Departamento de Matem\'aticas\\
Universidad Aut\'onoma de Madrid\\
28049-Madrid\\
Spain\\
\printead{e2}\\
\printead{u1}}
\end{aug}

\received{\smonth{8} \syear{2013}}
\revised{\smonth{11} \syear{2013}}

%
\begin{abstract}
A domain $S\subset{\mathbb R}^d$ is said to fulfill the \textit{Poincar\'e
cone property}
if any point in the boundary of $S$ is the vertex of a (finite) cone
which does not otherwise intersects the closure $\bar S$. For more than
a century, this condition has played a relevant role in the theory of
partial differential equations, as a shape assumption aimed to ensure
the existence of a solution for the classical Dirichlet problem on $S$.
In a completely different setting, this paper is devoted to analyze
some statistical applications of the Poincar\'e cone property (when
defined in a slightly stronger version). First, we show that this
condition can be seen as a sort of generalized convexity: while it is
considerably less restrictive than convexity, it still retains some
``convex flavour.'' In particular, when imposed to a probability
support $S$, this property allows the estimation of $S$ from a random
sample of points, using the ``hull principle'' much in the same way as
a convex support is estimated using the convex hull of the sample
points. The statistical properties of such hull estimator (consistency,
convergence rates, boundary estimation) are considered in detail.
Second, it is shown that the class of sets fulfilling the Poincar\'e property
is a $P$-Glivenko--Cantelli class for any absolutely continuous
distribution $P$ on~${\mathbb R}^d$. This has some independent interest
in the theory of empirical processes, since it extends the classical
analogous result, established for convex sets, to a much larger class.
Third, an algorithm to approximate the cone-convex hull of a finite
sample of points is proposed and some practical illustrations are given.
\end{abstract}

%


\pdfkeywords{62G05, 62G20, Poincare property, Glivenko-Cantelli classes, set estimation}

\end{frontmatter}

\section{Introduction}\label{sec:intro}

The Poincar\'e cone property (PCP) is a regularity condition for sets
in the Euclidean space. It has been used in mathematics (in partial
differential equations and Brownian motion theory) since more than a
century. We are concerned here with some new applications of this
property in statistics and probability. Let us begin by formally
establishing this condition, as well as other related notions we will use.

\subsection{The Poincar\'e property: Some history}

The standard version of PCP, which can be found in many books dealing
with potential theory
or Brownian motion is as follows; see, for example, M{\"o}rters and Peres \cite{mor10},
page 68:

\begin{definition}\label{def:pcp}
A domain $S\subset{\mathbb R}^d$
satisfies the \emph{Poincar\'e cone property} at $x\in\partial S$ if
there exist a cone $C(x)$ with vertex at $x$ and a number $h>0$ such that
%
\begin{equation}
C(x)\cap B(x,h)\subset S^c,\label{pz}
\end{equation}
where $B(x,h)$ denotes the open ball with center $x$ and radius $h$.
\end{definition}

The interest of this property is mainly associated with the so-called
Dirichlet problem which consist of finding a function $u$, harmonic on
$S$ (i.e., $\nabla^2 u=0$ on $S$) such that the restriction of $u$ to
the boundary $\partial S$ coincides with a given continuous function
$f$. This problem was posed by Gauss in 1840. During some years, it was
believed (from a conjecture due to Gauss himself) that the problem had
always a solution; however, this is not the case, unless some
regularity assumptions are imposed on $S$. In 1899, Poincar\'e showed
that a solution does exist whenever every point in $\partial S$ lies on
the surface of a sphere which does not otherwise intersects the closure
$\bar S$. In 1911, Zaremba showed that this ``outer sphere
condition'' proposed by
Poincar\'e could be weakened by replacing the sphere with a cone, as
indicated in Definition~\ref{def:pcp}. For this reason, condition
(\ref
{pz}) is sometimes also called Poincar\'e--Zaremba property
(e.g., Gilbarg and Trudinger \cite{gil77}) or just Zaremba's condition (Karatzas and Shreve \cite{kar88}, page~250).

Further details on the use of this classical property, its history and
its beautiful connections with the theory of Brownian motion can be
found, for example, in Kellogg \cite{kel29},
Gilbarg and Trudinger \cite{gil77}, Karatzas and Shreve \cite{kar88} and M{\"o}rters and Peres \cite{mor10}.

The intuitive meaning of Definition~\ref{def:pcp} is quite clear: if,
for a point $x\in\partial S$, we can always construct a ``finite
outside cone'' $C(x)\cap B(x,h)$ with vertex $x$, then we are typically ruling
out the existence of a sharp inward peak at $x$. A classical example of
a set
not fulfilling this condition is the so-called \textit{Lebesgue Thorn},
which is expressively described as follows in Kellogg \cite{kel29}, page 285:
\textit{Suppose we take a sphere with a deformable surface and at one of
its points push in a very sharp spine (\ldots)}. This set was first
proposed by Lebesgue in 1913 as a counterexample to show that the
Dirichlet problem is not always solvable.

\subsection{From the Poincar\'e property to cone-convexity: Our main
definitions}

As established in Definition~\ref{def:pcp}, the Poincar\'e cone
property is \textit{pointwise} in the sense that the opening angle
$\rho$ of the cone $C$ and the radius $h$ of the ball $B(x,h)$ in the
``$x$-cornet'' $C(x)\cap B(x,h)$ of condition (\ref{pz}) might depend
on $x$. For our statistical applications, we will need the condition
(\ref{pz}) to hold uniformly in $x$. Also, we will not be restricted to
assume that $S$ is a domain (i.e., an open connected set) since we are
interested in using the Poincar\'e property for support of probability
measures which, by definition, are closed sets.

So, in summary, the basic concepts we are going to handle arise as the
following strengthened versions of Definition~\ref{def:pcp}.

\begin{definition}\label{def:rcc}
We will say that the set $S\subset{\mathbb R}^d$ is \emph{$\rho
$-cone-convex} ($\rho$-cc), for some $\rho\in(0,\pi]$, if there
exists $h>0$ such that for all $x\in\partial S$ there is an open cone
$C_{\rho}(x)$ with opening angle $\rho$ and vertex $x$ such that the condition
%
\begin{equation}
C_\rho(x)\cap B(x,h)\subset S^c,\label{pzz}
\end{equation}
holds. When the above condition is satisfied for a specified $h>0$ we
will also say that
$S$ is \emph{$\rho,h$-cone-convex}.
\end{definition}

In informal terms, we could compare this definition with the standard
characterization of
convex sets in terms of supporting hyperplanes: if the (closed) set $S$
is convex, then for each $x\in\partial S$ there is a supporting
hyperplane $H=H(x)$ passing through $x$. Conversely, if $S$ is closed
with nonempty interior, the existence of a supporting hyperplane for
each $x\in\partial S$ implies that $S$ is convex. In Definition~\ref
{def:rcc}, we have replaced the supporting hyperplanes with
``supporting cones'' of
the form $C_{\rho,h}(x)=C_\rho(x)\cap B(x,h)$. Thus, for $\rho=\pi$ and
$h=\infty$, we would get as a particular case the hyperplane supporting
property for convex sets. Observe, however, that $\rho$-cone convexity
is a much more general condition than convexity as it allows the set
$S$ to have holes and inward peaks, as long as they are not too sharp
(the ``sharpness'' being limited by the angle $\rho$).

\subsection{Applications to set estimation}\label{aplicaciones}

We will explore here the applicability of Poincar\'e cone property from
a completely different point of view, mostly related with the problem
of \textit{set estimation}\, which basically deals with the reconstruction
of a set $S$ from a random sample points. See, for example, Cuevas and Fraiman \cite
{cue10} for an overview. Typically, in set estimation very little can
be said about the target set $S$ (beyond some simple results of
consistency) on the basis of the available sample information, unless
some relatively strong shape restrictions are imposed on $S$. Of course
such assumptions entail some loss in generality but, in return,
a~wealth of valuable results (estimation of the boundary and the boundary
measure, rates of convergence, etc.) are typically obtained.

\textit{The convex case and the ``hull mechanism.''}
The use we will make of the Poincar\'e condition, via $\rho
$-cone-convexity, is better explained from the perspective of other
more popular related properties.
The most obvious one is convexity. If the support $S$ is assumed to be
convex, then the natural estimator of $S$ from a sample $X_1,\ldots
,X_n$ is the convex hull $S_n(X_1,\ldots,X_n)$, that is the minimal
convex set including the sample points. The study of the convex hull
itself (even not considering its properties as an estimator of~$S$)
seems to be an inexhaustible subject of research in geometric
probability: increasingly sophisticated results on the distributions of
several random variables (number of vertices, area, perimeter,
probability content, number of sides, etc.) associated with $S_n$ have
been considered in the last fifty years. The survey paper by Reitzner \cite
{rei10} provides an excellent up-to-date account of these topics.

The properties of the convex hull $S_n$ as an estimator of the support
$S$, and in particular the convergence rates for the Hausdorff distance
$d_H(S_n,S)$, are studied by D{\"u}mbgen and Walther
\cite{dum96} among others.

Thus, convexity is the prototypical example where the \textit{hull
mechanism} (that is to define our estimator as the ``minimal one
including the sample and fulfilling a desired shape property'') can be
successfully used. It is natural to ask whether in other cases, under
more general assumptions than convexity, the hull mechanism could also work.

\textit{From $r$-convexity to $\rho$-cone convexity}.
The so-called $r$-\textit{convexity property} provides an interesting
example: a closed set $S$ is said to be $r$-\textit{convex} if it can
be expressed as the intersection of a family of complements of balls
with radius $r>0$. More precisely, $S$ is $r$-convex if and only if
%
\begin{equation}
S=\bigcap_{\{y: B(y,r)\cap S=\varnothing\}}B(y,r)^c.\label{rconv}
\end{equation}

It is easy to check that any convex set is also $r$-convex for all
$r>0$ but, clearly, $r$-convexity is a much milder restriction. In
particular, it allows for smooth or ``round gulfs,'' and even holes, in
the set.

The study of this property dates back to Perkal \cite{per56}; see also
Walther \cite
{wal97} for further statistical insights on this concept. From a
statistical viewpoint, the interesting fact is that, if a set $S$ is
assumed to be $r$-convex, then it can be (asymptotically) recovered
from a random sample by just considering the $r$-\textit{convex hull}
of the data points as a natural estimator.

The effective calculation of this $r$-convex hull is much more involved
than that of the ordinary convex hull. The R-package \textit{alphahull}
provides a practical implementation for the case $d=2$; see Pateiro-L{\'o}pez and Rodr{\'{\i
}}guez-Casal \cite
{pat10}. Whereas the ordinary convex hull of a sample in the plane is
always a polygon, the boundary of the $r$-convex hull is made of arcs
of $r$-circumferences plus, perhaps, some isolated points;
see Figure~\ref{fig4}
in Section~\ref{numerical-sec} for an example. More information on
statistical properties, examples and applications of the $r$-convex
hull can be found in Rodr{\'{\i}}guez-Casal \cite{rod07}, Pateiro-L{\'o}pez and Rodr{\'{\i
}}guez-Casal \cite{pat08} and Berrendero, Cuevas and Pateiro-L{\'
o}pez \cite{ber12}.
Cuevas, Fraiman and Pateiro-L{\'
o}pez \cite{cue12} provide further results on the estimation of an
$r$-convex support as well as also some insights regarding the
comparison of $r$-convexity with other better known properties such as
positive reach (Federer \cite{fed59}) and the above mentioned (uniform) outer
$r$-sphere property. In particular, $r$-convexity is shown to be
slightly stronger than the ``rolling'' outer $r$-sphere property: for
every point in the boundary of $S$ there exists a ball touching that
point whose interior is included in $S^c$.

\begin{figure}

\includegraphics{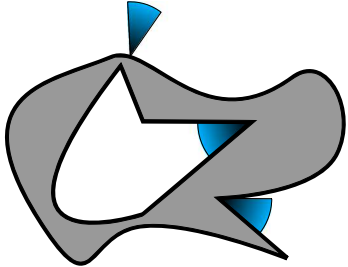}

\caption{A general $\rho$-cone-convex set with inward peaks and holes.}
\label{light}
\end{figure}

In this paper, we replace the outer balls by ``outer cones'' (in the
spirit of Poincar\'e's definition). This led us in a natural way to the
cone-convexity notion introduced in Definition~\ref{def:rcc}. In a
similar vein, expression (\ref{rconv}) suggests the following
cone-based analogue notion.
%
\begin{definition}\label{rccc}
A closed set $S\subset\mathbb{R}^d$ is said to be \emph
{cone-convex by
complement} with parameters $\rho\in(0,\pi]$, $h>0$ ($\rho,h$-ccc) if
and only if
%
\begin{equation}
S=\bigcap_{\{y: C_{\rho,h}(y)\cap S=\varnothing\}} \bigl(C_{\rho
,h}(y)
\bigr)^c,\label{rhoccc}
\end{equation}
where $C_{\rho,h}(y)$ denotes a finite cone with vertex $y$, of type
$C_{\rho,h}(y)=C_\rho(y)\cap B(y,h)$.
\end{definition}

In informal terms, one could say that a closed set $S$ is convex by
complement (with parameters $\rho$, $h$) if any point $x\notin S$ can
be separated from $S$ by a finite cone
$C_{\rho,h}(y)$, with opening angle $\rho$ and height $h$, which
contains $x$.

Thus, in summary, the $\rho$-cone-convexity properties considered in
this paper are two generalizations of the notion of $r$-convexity where
the balls are replaced by finite cones; see Figure~\ref{light}. These
generalizations allow us to consider\vadjust{\goodbreak} much more general sets with
rougher boundaries. To be more precise, the cone convexity by
complement is a direct extension of the notion of
$r$-convexity, by replacing the balls of radius $r$ with the $\rho
,h$-cones. Likewise, the $\rho,h$-cone-convexity is a generalization of
the ``outer rolling property'' commented above (i.e., any boundary
point of $S$ has a touching ball whose interior is included in $S^c$;
see Cuevas, Fraiman and Pateiro-L{\'
o}pez \cite{cue12} for details). However, whereas the $r$-convexity
implies the outer rolling ball property (see Proposition~2 in Cuevas, Fraiman and Pateiro-L{\'
o}pez \cite
{cue12}), the analogous implication does not hold for the cone-convex
case: see Proposition~\ref{inclusiones} below.

We will show that, in spite of this gain in generality, the ``hull
principle'' still works for the $\rho$-cone-convex properties, so that
it can be also employed for estimation purposes. This means that a
$\rho
$-cc (or a $\rho$-ccc) support
$S$ can be estimated, from a random sample drawn on $S$, just using the
corresponding $\rho$-cone-convex hull of the sample points.

In addition, a relevant property (see Theorem~\ref{GCtheorem}) is also
shown for the class of $\rho$-cone convex sets: whereas this class is
considerably broader than that of convex sets, it is still a
Glivenko--Cantelli class. This represents a generalization of the recent
similar result proved by Cuevas, Fraiman and Pateiro-L{\'
o}pez \cite{cue12} for the case of $r$-convex sets.
An application is given in Theorem~\ref{GC-consistencia}.

\subsection{Some notation. The organization of this paper}\label{notation}

With some notational abuse, a ``cornet'' obtained by intersecting an
infinite cone with a ball centered at its vertex is called itself as
\emph{a (finite) cone}.
A set of this type is thus defined by the vertex $x$, a unit vector
$\xi
$ indicating the axis of the cone, an angle $\rho\in(0,\pi]$ indicating
the opening angle and a positive number $h>0$ corresponding to the
radius of the intersecting ball.

Thus, in precise terms, an infinite cone is defined by
\[
C_{\rho}(x)= \biggl\{z\in{\mathbb R}^d, z\neq x\dvtx \biggl
\langle\xi, \frac{z-x}{\|z-x\|} \biggr\rangle>\cos(\rho/2) \biggr\},
\]
and, for $h>0$, we will denote
$C_{\rho,\xi,h}(x)=B(x,h)\cap C_{\rho}(x)$.
The subindices, especially $\xi$, will be omitted when convenient, so
the notation $C_{\rho,h}(x)$ is often used for finite cones.

The class of nonempty compact sets $S\subset{\mathbb R}^d$ satisfying,
for a given $h>0$, the $\rho$-cone-convex
condition (\ref{pzz}) established in Definition~\ref{def:rcc} will be
denoted by
${\mathcal C}_{\rho,h}$.

Also, the class of nonempty compact sets $S\subset{\mathbb R}^d$
satisfying, for a given $h>0$, the $\rho$-ccc
condition established in Definition~\ref{rccc} will be denoted by
$\tilde{\mathcal C}_{\rho,h}$.
If $x\in\R{d}$ and $S\subset\R{d}$, $S\neq\varnothing$, the distance
from $x$ to $S$ is $d(x,S)=\inf\{\Vert x-s\Vert\dvtx  s\in S\}$. The
Lebesgue measure on ${\mathbb R}^d$ will be denoted by $\mu$.
Given a bounded set $A\subset{\mathbb R}^d$ and $\varepsilon>0$,
$B(A,\varepsilon)$ will denote the parallel set
$B(A,\varepsilon)=\{x\in{\mathbb R}^d\dvtx d(x,A)\leq\varepsilon\}$.
Note that, according to this notation, $B(\{x\},\varepsilon)$ coincides
with the closed ball
centered at $x$ with radius $\varepsilon$, not with the open ball
$B(x,\varepsilon)$.

Given two compact nonempty sets $A,B\subset{\mathbb R}^d$, the
\textit{Hausdorff distance} or
\textit{Hausdorff--Pompein distance} between
$A$ and $C$ is defined by
%
\begin{equation}
d_H(A,C)=\inf \bigl\{\varepsilon>0\dvtx \mbox{such that } A\subset
B(C, \varepsilon) \mbox{ and } C\subset B(A,\varepsilon) \bigr\}.\label{Hausdorff}
\end{equation}
The class $\mathcal{M}$ of compact nonempty sets of ${\mathbb R}^d$,
endowed with the distance $d_H$
is known to be a complete separable metric space; see, for example,
Rockafellar and Wets \cite{roc09}, Chapter~4. Moreover, any class of uniformly bounded
subsets in such space is relatively compact with respect to $d_H$. So,
any bounded sequence of compact nonempty subsets of ${\mathbb R}^d$
has a convergent subsequence.
For a given Borel measure $\nu$, define also the pseudometric
$d_\nu(A,C)=\nu(A\setminus C)+\nu(C\setminus A)$.

The rest of this paper is organized as follows. In Section~\ref{hulls}, we analyze the notions of ``convex hulls'' associated
with the concepts of cone-convexity introduced above. Some general
properties of convergence for sequences of cone-convex sets are
obtained in Section~\ref{convergencia}. Section~\ref{GC} is
devoted to show that the class of cone-convex sets is a
Glivenko--Cantelli class. This has some independent interest in the
theory of empirical processes, since it extends the classical
analogous result, established for convex sets, to a much larger
class. The estimation (consistency and convergence rates) of
cone-convex sets using the corresponding cone-convex hull of the
sample is considered in Section~\ref{estimacion}. A stochastic
algorithm to approximate the cone-convex hull by complement of a
sample is provided in Section~\ref{algoritmo}. The behavior of
this algorithm is illustrated with some examples and simulations
in Section~\ref{numerical-sec}. Some final comments and
suggestions for further work are given in Section~\ref{final-sec}.

\section{The notion of cone-convex hull}\label{hulls}

We now define the concept of cone-convex hull corresponding to the
notion we have introduced of cone-convexity. In fact, we will need to
distinguish between the ``cone
convex hull'' and the ``cone convex hull by complement'' which, unlike the
classical convex case, do not coincide in general for cone-convexity.
Let $S\subset{\mathbb R}^d$ be a bounded set.

\begin{definition} \label{envolvente cono convexa}
(a) The \emph{$\rho,h$-cone-convex hull} ($\rho,h$-cc) of $S$ or, just, the
cone-convex hull of $S$, is defined by
%
\begin{equation}
\label{ii} {\mathbb C}_{\rho,h}(S)=\bigcap
_{S\subset B_t, B_t \in\C_{\rho,h}}B_t.
\end{equation}

\begin{longlist}[(b)]
\item[(b)]The \emph{$\rho,h$-cone-convex hull by complement} ($\rho
,h$-ccc) of $S$, is defined as the intersection of the complements of
those (open, finite) cones $C_{\rho,\xi,h}$ which do not intersect $S$.
We will denote it by $\tilde{\mathbb C}_{\rho,h}(S)$.
Note that $\tilde{\mathbb C}_{\rho,h}(S)$ can be also expressed as
%
\begin{equation}
\label{iii} \tilde{\mathbb C}_{\rho,h}(S)=\bigcap
_{S\subset B_t, B_t \in\tilde\C
_{\rho,h}}B_t.
\end{equation}
\end{longlist}
\end{definition}

\begin{proposition}\label{inclusiones}
Given $\rho\in(0,\pi]$ and $h>0$, we have the following relations
between the two classes $\C_{\rho,h}$ and $\tilde\C_{\rho,h}$ of
cone-convex sets, introduced in Section~\ref{notation}, and the
corresponding cone-convex hulls ${\mathbb C}_{\rho,h}(S)$ and $\tilde
{\mathbb C}_{\rho,h}(S)$ defined for any bounded $S\subset{\mathbb R}^d$.
\begin{longlist}[(a)]
\item[(a)] Neither of the classes $\C_{\rho,h}$ and $\tilde\C
_{\rho,h}$
is included in the other.

\item[(b)] ${\mathbb C}_{\rho,h}(S)\in\C_{\rho,h}$ and
$\tilde{\mathbb C}_{\rho,h}(S)\in\tilde\C_{\rho,h}$.

\item[(c)] If $\rho^\prime\leq\rho$ and $h^\prime\leq h$, then
$\C
_{\rho,h}\subset\C_{\rho^\prime,h^\prime}$ and $\tilde\C_{\rho
,h}\subset\tilde\C_{\rho^\prime,h^\prime}$. Also,\break  ${\mathbb
C}_{\rho
^\prime,h^\prime}(S)=S$ for all $S\in\C_{\rho,h}$ and $\tilde
{\mathbb
C}_{\rho^\prime,h^\prime}(S)=S$, for all $S\in\tilde\C_{\rho,h}$.

\item[(d)] Let $S\subset\R{d}$ be bounded. For all $\rho\in(0,\pi]$
and $h>0$, let us define $h'=h/2 \sin(\rho/2)$, $\rho'= (\pi-\rho)/2$
if $\rho>\pi/3$ and $\rho'=\rho$ if $\rho\leq\pi/3$. Then
$\tilde
{\mathbb{C}}_{\rho,h} (S)\in\C_{\rho',h'}$. As a consequence, if $S$
is $\rho, h$-ccc, then it is also
$\rho', h'$-cc.
\end{longlist}
\end{proposition}
\begin{pf}
(a) Let $S$ denote the closed unit ball in ${\mathbb R}^2$ and let $C$
be any (open) cone with vertex at the origin $(0,0)$, opening angle
$\rho\in(0,\pi/2]$ and height $h<1/2$. Then it is readily seen that
the set
$S\setminus C$ belongs to $\tilde\C_{\rho,h}$ but not to $\C_{\rho,h}$
since condition (\ref{pzz}) fails for all $x\in\partial C$, $x\neq(0,0)$.

\begin{figure}

\includegraphics{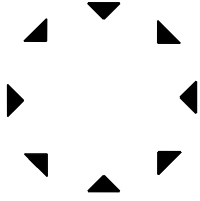}

\caption{An example for which the envelopes
${\mathbb C}_{\pi/4,3}(S)$ and $\tilde{\mathbb C}_{\pi/4,3}(S)$ do not
coincide.}\label{triang}
\end{figure}

Also, the set $E=E_1\cup E_2\subset{\mathbb R}^2$ where $E_1$ is the
graph of the function $f(t)=t^2$ on $t\in[0,1]$ and $E_2=\{(t,0)\dvtx 0\leq t\leq1\}$ belongs to $\C_{\rho,h}$ for all $\rho\in(0,\pi/2]$,
$h>0$ but $E\notin\tilde\C_{\rho,h}$ for any $\rho$.

There are also counterexamples of sets $S$ with nonempty interior such
that $S\in\C_{\rho,h}\setminus\tilde\C_{\rho,h}$: let $S$ be
union of
the triangle $T$ with vertices
$(1,0)$, $(1+s,s)$ and $(1+s,-s)$, where $s=(\sqrt{2}-1)/2$, plus the
seven congruent triangles obtained from $T$ by applying a rotation
around $(0,0)$ with angle $\pi/4$; see Figure~\ref{triang}.

This set is $\pi/4, h$-cone-convex for any $h$. However, the origin
cannot be separated from $S$ by any cone $ C_{\pi/4, 3}(y)$ since any
such cone should contain at least one of the vertices of the congruent
triangles of $S$.

(b) By definition, we have
\[
{\mathbb C}_{\rho,h}(S)=\bigcap_{S\subset B_t, B_t \in\C_{\rho,h}}
B_t.
\]
Given $x\in\partial{\mathbb C}_{\rho,h}(S)$ we want to find a
$C_{\rho
,h}(x)$ with $C_{\rho,h}(x)\subset{\mathbb C}_{\rho,h}(S)^c$.
We have $B(x,1/n) \cap B_t \neq\varnothing$ for all $B_t$ such that
$S\subset B_t$ and $B_t \in\C_{\rho,h}$. Moreover, $B(x,1/n)\cap B_n^c
\neq\varnothing$ for some $B_n$ with $S\subset B_n$ and $B_n\in\C
_{\rho
,h}$. Given $z_n\in B(x,1/n)\cap\partial B_n$, we have $z_n\rightarrow
x$. Since $z_n\in\partial B_n$ and $B_n \in\C_{\rho,h}$, there must
exist a cone $C_{\rho,\xi_n,h}(z_n)\subset B_n^c$. Since $\|\xi_n\|=1$
there exists $\xi=\lim_n \xi_n$ (for some subsequence $\xi_n$) and also
$z_n\rightarrow x$. We thus have that $\overline{C_{\rho,\xi
_n,h}(z_n)}$ converges in the Hausdorff metric to $\overline{C_{\rho
,\xi
,h}(x)}$.
We only must check $C_{\rho,\xi,h}(x)\cap S=\varnothing$. Indeed,
otherwise we would have some $s\in S\cap C_{\rho,\xi,h}(x)$ with
$ \langle\frac{s-x}{\|s-x\|},\xi\rangle=\cos(\rho/2)+\delta$,
for some $\delta>0$. Then
$ \langle\frac{s-z_n}{\|s-z_n\|},\xi_n \rangle>\cos(\rho/2)$,
for $n$ large enough, which contradicts $C_{\rho,\xi,h}(z_n)\cap
S=\varnothing$.

The second statement follows directly from the expression (\ref{iii}).

(c) This is a direct consequence of the definitions of the classes and
the respective hulls.

(d) We want to find $\rho',h'$ such that for any $z\in\partial\tilde
{\mathbb{C}}_{\rho,h} (S)$ we have
a cone $C_{\rho',h'}(z)$ with $C_{\rho',h'}(z)\subset\tilde{\mathbb
{C}}_{\rho,h} (S)^c$. Now, take a sequence $z_n\to z$ with $z_n\in
\tilde{\mathbb{C}}_{\rho,h} (S)^c$. Using the definition of
$\tilde{\mathbb{C}}_{\rho,h} (S)$, there is a sequence of cones
$C_{\rho,h}(u_n)$, disjoint with $\tilde{\mathbb{C}}_{\rho,h} (S)$,
such that $z_n\in C_{\rho,h}(u_n)$. Take a further subsequence such
that the cones $C_{\rho,h}(u_n)$ are convergent, that is, $C_{\rho,h}(u_n)
\to C_{\rho,h}(u)$.
By construction, we have that $z\in\partial C_{\rho,h}(u)$ and
$C_{\rho
,h}(u)\cap\tilde{\mathbb{C}}_{\rho,h} (S)=\varnothing$. It only remains
to show that
$C_{\rho',h'}(z)\subset C_{\rho,h}(u)$.
This is readily seen for the given values of $\rho'$ and $h'$; see
Figure~\ref{Prop1d}. Indeed, the result is immediate for any $z$ at a
distance $h/2$ from the vertex $x$. For the remaining $z\in\partial
C_{\rho,h}(z)$ simple translations of this cone, perhaps combined with
a rotation of the cone axis provide the required $C_{\rho',h'}(z)$.
\end{pf}

\begin{figure}

\includegraphics{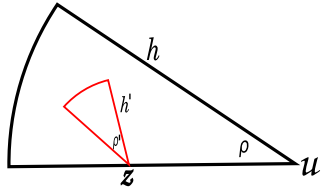}

\caption{For any $z\in\partial C_{\rho,h}(u)$ there is a cone
$C_{\rho
',h'}(z)\subset C_{\rho,h}(u)$.}
\label{Prop1d}
\end{figure}

\textit{Practical consequences in estimation problems}.
As a conclusion of the above result, we have two ways of estimating a
cone-convex set $S$ from a sample $\aleph_n=\{X_1,\ldots,X_n\}$ drawn
from a distribution whose support is $S$. If we assume that
$S$ is $\rho,h$-cone-convex, then the natural estimator of $S$ would be
${\mathbb C}_{\rho,h}(\aleph_n)$.
When $S$ is assumed to be $\rho,h$-cone-convex by complement, then
$\tilde{\mathbb C}_{\rho,h}(\aleph_n)$ would be the natural estimator
of $S$.

The difference between both notions of cone-convexity is mainly
technical. In fact, there is a considerable overlapping between the
classes $\tilde\C_{\rho,h}$ and $\C_{\rho,h}$: most sets found in
practice fulfilling one of these conditions will also satisfy the other one.
For example, as pointed out above, if $S$ is $r$-convex (i.e., it can
be expressed as the intersection of the complements of a family of
$r$-balls), then $S$ fulfils both Definitions \ref{def:rcc} and \ref
{rccc} of $\rho,r$-cone-convexity. In those cases, both envelopes
${\mathbb C}_{\rho,h}(\aleph_n)$ and $\tilde{\mathbb C}_{\rho
,h}(\aleph
_n)$ can be used.

We will analyze the asymptotic properties of both estimators, but the
envelope $\tilde{\mathbb C}_{\rho,h}(\aleph_n)$
is easier to approximate via an stochastic algorithm; see Section~\ref{algoritmo}.

In practice, the correct choice of the parameters $\rho$, $h$ will
depend on prior assumptions on the nature of the sets under study.
Note, however, that result (c) in Proposition~\ref{inclusiones}
guarantees that an exact knowledge of the ``optimal'' (maximal) values
of these parameters is not needed, in the sense that a conservative
(small) choice of $\rho$ and $h$ would do the job.

\subsection{Lighthouses}\label{lighthouses}

A particular case of Definition~\ref{def:rcc} deserves attention as it
represents a much more direct extension of the convexity notion: if
condition (\ref{pzz})
holds for all $h>0$ then, for each point in $\partial S$ we can find an
infinite supporting cone on $x$. In this case, for $\rho=\pi$ condition
(\ref{pzz}) amounts to the supporting hyperplane condition.
The formal definition would be as follows.

\begin{definition}
We will say that $S$ is a \emph{$\rho$-lighthouse set} when condition
(\ref{pzz}) holds for all $h>0$ or, equivalently, when for all $x\in
\partial S$ there is an open cone $C_\rho(x)$ based on $x$ with opening
angle $\rho$ such that
%
\begin{equation}
C_\rho(x)\subset S^c.\label{pzzz}
\end{equation}
\end{definition}

It is clear that the class ${\mathcal C}_{\rho}$ of compact sets in
${\mathbb R}^d$ fulfilling condition (\ref{pzzz}) is much broader than
the class ${\mathcal C}$ of compact convex sets but, definitely, much
smaller than any family ${\mathcal C}_{\rho,h}$ of $\rho,h$-cone-convex
sets since condition (\ref{pzzz}) would typically exclude the presence
of holes in $S$ [provided that $S=\overline{\operatorname{int}(S)}$]. In
graphical terms, (\ref{pzzz}) imposes the possibility of illuminating
the space around $S$ with a full beam of light from any boundary point
$x$ in $S$. This accounts for the term ``lighthouse.''

Likewise, the cone-convex hull notions introduced in Definition~\ref
{envolvente cono convexa} can be readily adapted to the lighthouse
sets just replacing the finite cones in both (\ref{ii}) and (\ref{iii})
by infinite (unbounded) cones.

In the next three sections, we will focus on the general case of $\rho
,h$-cone-convex sets (for a finite $h$) but
our results might be translated to the case of $\rho$-lighthouses.
Apart from the simplicity and intuitive appeal of the lighthouse
condition, the inference on the parameter $\rho$ is expected to be much
easier in this case. However, this topic is not considered here.

\section{Convergence properties}\label{convergencia}

Let us start with a simple regularity property of cone-convex sets.

\begin{proposition} \label{muS} If $S \in{\mathcal C}_{\rho,h}$, then
$\mu(\partial S)=0$.
\end{proposition}

\begin{pf} First note that $\partial S$ is a Borel set since $S$ is closed.
Now let us recall that a point $x\in S$ is said to have \textit{metric
density} 1 (see, e.g., Erd{\"o}s \cite{erd45}) if for all $\varepsilon>0$
there is some $\delta>0$ such that
$\mu(S\cap B(x,r) )>(1-\varepsilon)\mu(B(x,r) )$ for all
$r<\delta$. From Corollary~2.9 in Morgan \cite{mor00}, every set with
positive (Lebesgue) measure has at least a point with metric density 1.
This implies that we must have $\mu(\partial S)=0$ since
for all $x\in\partial S$ there exists an open cone $C_{\rho
,h}(x)\subset S^c$. Therefore, for all $r<h$, we have some $k_{\rho,h}<1$ such that
\begin{eqnarray*}
\mu \bigl( \partial S\cap B(x,r) \bigr)&\leq&\mu \bigl(S\cap B(x,r) \bigr)= \mu
\bigl( S\cap C_{\rho,h}(x)^c\cap B(x,r) \bigr)\\
&\leq&
k_{\rho,h} \mu \bigl(B(x,r) \bigr).
\end{eqnarray*}
\upqed\end{pf}

We now establish that the convergence of a sequence of $\rho
$-cone-convex sets entails the convergence of their respective
boundaries. This is an important regularity property. It essentially
says that we cannot have $\rho$-cone convex sets very close together if
the respective boundaries are far away from each other. A similar
property has been recently proved for sets fulfilling the rolling
condition [i.e., property (\ref{pz})] where $C(x)\cap B(x,h)$ is
replaced by an open ball $B(x,r)$ of a given radius $r$, whose boundary
contains $x$ (see
Theorem~3(a) in Cuevas, Fraiman and Pateiro-L{\'
o}pez \cite{cue12}; see also Ba{\'{\i}}llo and Cuevas \cite{bai01} for related
results for the case of star-shaped sets). Our Theorem~\ref{teo1} below
can be seen as a considerable extension of this result (since $\rho,h
$-cone convexity is a much less restrictive than the rolling property).

\begin{theorem} \label{teo1} Let $\{S_n\}\subset\C_{\rho,h}$ (or $\{
S_n\}\subset\tilde{\C}_{\rho,h}$) and let $S\subset{\mathbb R}^d$
be a
compact set such that $d_H(S_n,S)\rightarrow0$. Then
$d_H(\partial S_n, \partial S)\rightarrow0$.
\end{theorem}
\begin{pf} By contradiction, let us assume that $d_H(\partial S_n,
\partial S)\nrightarrow0$. Then we should have either (i) or (ii):
\begin{longlist}[(ii)]
\item[(i)] There exists $\eps>0$ such that for some subsequence $x_n
\in\partial S$ we have $d(x_n,\partial S_n)>\eps$.
\item[(ii)] There exists $\eps>0$ such that for some subsequence $x_n
\in\partial S_n$ we have $d(x_n,\partial S)>\eps$.
\end{longlist}

\indent Suppose that we have a sequence $x_n$ fulfilling (i). Since $S$
is compact, there exists a convergent subsequence, denoted again $x_n$.
Let $x$ be the limit of such subsequence.
Since $d_H(S_n,S)\rightarrow0$, we have\vadjust{\goodbreak} that, for $n$ large enough,
$d(x,S_n)\leq\eps/2$. Moreover, taking the infimum
on $y\in\partial S_n$ in $\|x_n-y\|\leq\|x_n-x\|+ \|x-y\|$, we have
that, eventually,
\mbox{$d(x, \partial S_n)\geq d(x_n,\partial S_n)\!-\!\|x\!-\!x_n\|>\eps/2$}.

On the other hand, since $d(x,S_n)\leq\eps/2$ and $d(x, \partial
S_n)>\eps/2$ we have $x\in \operatorname{int}(S_n)$ and $B(x,\eps/2)\cap
S_n^c=\varnothing$. Since $x\in\partial S$ and $S$ is compact, we may
take $y\in S^c$ such that $\|x-y\|<\eps/2$ and $d(y,S)>0$. But
$d_H(S_n,S)\rightarrow0$ and $d(y,S)>0$ entail $d(y,S_n)>0$
eventually, so $y\in B(x,\eps/2)\cap S_n^c$, in contradiction with
$B(x,\eps/2)\cap S_n^c=\varnothing$.

Let us now assume that we have (ii). Since $d_H(S_n,S)\to0$, we must
have $d(x_n,S)<\eps$ eventually which, together with $d(x_n,\partial
S)>\eps$, yields [by a similar reasoning to that in (i)] $x_n\in
\operatorname{int}(S)$ eventually.
Take now a convergent subsequence of $x_n$ (denoted again $x_n$) with
$x_n\rightarrow x\in S$ and $d(x,\partial S)>\eps/2$, that is, $B(x,\eps
/2)\subset S$.
Since $x_n\in\partial S_n$ and $S_n\in\C_{\rho,h}$, there exist
finite cones $C_{\rho,h}(x_n)$ with $C_{\rho,h}(x_n)\cap
S_n=\varnothing
$. Take $c_n\in S_n^c$, in the axis of $C_{\rho,h}(x_n)$, with $0<\|
x_n-c_n\|=\min\{h/2,\eps/4\}=k$. We may assume (taking, if necessary, a
further suitable subsequence) $c_n\rightarrow c$ and $\Vert c-x\Vert
\leq\eps/4$. Note that, by construction, $d(c_n,S_n)\geq k\sin(\rho
/2)$. Hence, $d(c,S_n)>k'$ eventually, for some constant $k'>0$. Since
$d_H(S_n,S)\to0$, we must have $c\in S^c$, in contradiction with $c\in
B(x,\eps/2)\subset S$.

For the case $\{S_n\}\subset\tilde{\C}_{\rho,h}$, the result follows
as a consequence of the previous case together with Proposition~\ref
{inclusiones}(d).
\end{pf}

The following result shows that the class $\C_{\rho,h}$ is
topologically closed.

\begin{theorem} \label{teo2}
Let $\{S_n\}\subset\C_{\rho,h}$ and let $S\subset{\mathbb R}^d$ be a
compact set such that $d_H(S_n,S)\rightarrow0$. Then,
$S\in\C_{\rho,h}$.
\end{theorem}
\begin{pf} Given $x\in\partial S$, we want to find a finite cone
$C_{\rho,h}(x)\subset S^c$. From Theorem~\ref{teo1}, we know that
$d_H(\partial S_n,\partial S)\rightarrow0$. Therefore, there is a
sequence $x_n\in\partial S_n$
such that $x_n\rightarrow x$.

Now the reasoning to find a cone $C_{\rho,h}(x)$ is similar to that of
Proposition~\ref{inclusiones}. By considering, if necessary, a suitable subsequence
of $x_n$, we
may find again a sequence of finite closed cones $\{\overline{C_{\rho
,\xi_n,h}(x_n)}\}$,
with $C_{\rho,\xi_n,h}(x_n)\subset S_n^c$, converging in the Hausdorff
metric to the closed cone
$\overline{C_{\rho,\xi,h}(x)}$ for some direction $\xi$ with \mbox{$\|\xi
\|
=1$}. We now check that this is the cone we are looking for, that is,
$C_{\rho,\xi,h}(x)\cap S=\varnothing$.

Suppose, by contradiction, that there is $z\in C_{\rho,\xi,h}(x)\cap
S$. Since $\xi_n\rightarrow\xi$ and $x_n\rightarrow x$, we may take
$\eps>0$ such that for $n$ large enough $B(z,\eps)\subset C_{\rho
,\xi
,h}(x_n)$. This entails $d(z,S_n)\geq\eps$ infinitely often, which
contradicts $d_H(S_n,S)\to0$.
\end{pf}

\section{The Glivenko--Cantelli property and the Poincar\'e
condition}\label{GC}

Let $X_1,\break  \ldots, X_n,\ldots$ be a sequence of i.i.d., ${\mathbb
R}^d$-valued random variables defined on a probability space $(\Omega,
\F, \pr)$. Denote by ${P}$\vadjust{\goodbreak} the common distribution of the $X_i$'s
on ${\mathbb R}^d$ and by $\pr_n$ the empirical distribution associated
with the first $n$ sample observations $X_1,\ldots,X_n$.

The main result of this section will be established in Section~\ref{GCresult} below.
In order to view this result from an appropriate
perspective, we next summarize some basic facts on the
Glivenko--Cantelli (GC) property.

\subsection{The Glivenko--Cantelli property: Some background}\label{GCsum}

A class $\A$ of Borel subsets of ${\mathbb R}^d$ is said to be a
\textit{$P$-Glivenko--Cantelli class} whenever
%
\begin{equation}
\sup_{A\in\A}\bigl | \pr_n(A)- {P}(A) \bigr|\rightarrow0\qquad
\mbox{a.s.}\label{GC-def}
\end{equation}
These classes are named after the classical Glivenko--Cantelli theorem
that establishes (\ref{GC-def}) for
the case $d=1$ when $\A$ is the class of half-lines of type $(-\infty
,x]$, with $x\in{\mathbb R}$.

In informal terms, a class $\A$ of sets is a GC-class if it is small
enough as to ensure the uniform validity of the strong law of large
numbers on $\A$. The study of GC-classes is a classical topic in the
theory of empirical processes. See, for example, Shorack and Wellner \cite{sho86} and
van~der Vaart \cite{vaa98}, Chapter~19, for detailed accounts of this theory and its
statistical applications. For example, Theorems 19.4 and 19.13 in
van~der Vaart \cite{vaa98} provide sufficient conditions for a class being
$P$-Glivenko--Cantelli. These conditions are expressed in terms of
entropy conditions which, in some sense, quantify the ``size'' of the
class $\A$. Chapters~12 and 13 in the book Devroye, Gy{\"o}rfi and
Lugosi \cite{dev96} provide an
insightful presentation of the Vapnik--Cervonenkis approach to the
study of GC-classes. In that approach, the GC-property is obtained
through an exponential bound for the probability ${\mathbb P}\{\sup_{A\in\A} |\pr_n(A)-P(A) |>\varepsilon\}$. If the series (in $n$)
of these bounds is convergent, then the classical Borel--Cantelli
lemma, leads to (\ref{GC-def}). However, this approach fails for some
important classes $\A$ as it requires the finiteness of the so-called
\textit{Vapnik--Cervonenkis (VC)} dimension of $\A$ (see Definition~12.2 in Devroye, Gy{\"o}rfi and
Lugosi \cite{dev96}). A different approach, which can be used in fact
for the study of GC-classes of functions is given by Talagrand \cite{tal87}.
That approach can be used to establish the GC-property in some
situations where the VC-dimension of $\A$ is infinite. Thus, it can be
proved (as a consequence of Theorem~5 in Talagrand~\cite{tal87}) that, given a
probability $P$ in ${\mathbb R}^d$, the class ${\mathcal C}$ of closed
convex sets in ${\mathbb R}^d$ such that $P(\partial C)=0$ is a $P$-GC-class.

We shall use here a different, older approach to the GC-problem due to
Billingsley and Tops{\o}e \cite{bil67}. Billingsley--Tops{\o}e approach can be used to prove a
property called $P$-uniformity, which is in fact more general than the
Glivenko--Cantelli condition, as it applies to general sequences of
probability measures (not necessarily empirical measures). A class of
sets $\A$ is said to be a $P$-uniformity class if
%
\begin{equation}
\sup_{A\in\A} \bigl|P_n(A)-P(A) \bigr|\to0 \label{Puniform}
\end{equation}
holds for every sequence $P_n$ of probability measures converging
weakly to $P$ (this is denoted $P_n\stackrel{w}{\rightarrow}P$) in
the sense
that $P_n(B)\to P(B)$ for every Borel set $B$ such that $P(\partial
B)=0$ [which of course happens, a.s. for $P_n=\pr_n$; therefore, (\ref
{Puniform}) implies (\ref{GC-def})].

We will use a result in Billingsley and Tops{\o}e \cite{bil67}, Theorem~4, according to
which if
$\mathcal{A}$ is a $P$-continuity class of Borel sets in ${\mathbb
R}^d$ [i.e., $P(\partial A)=0$ for every $A\in\A$] then the compactness
of the class $\{\partial A\dvtx A\in\A\}$, in the Hausdorff topology, is
a sufficient condition for
$\mathcal{A}$ to be a $P$-uniformity class.

In Theorem~5 of Cuevas, Fraiman and Pateiro-L{\'
o}pez \cite{cue12}, it is proved, using this result, that
any class $\A$ of nonempty closed sets $A\subset{\mathbb R}^d$,
uniformly bounded (i.e., all of them included in some common compact
$K$) and fulfilling $\operatorname{reach}(A)\geq r$, for some given $r>0$, is a
$P$-uniformity class whenever $P$ is a absolutely continuous with
respect to the Lebesgue measure. Here, $\operatorname{reach}(A)$ denotes the
supremum (possibly) of those values $s$ such that any point $x$ whose
distance to $A$ is smaller than $s$ has just one closest point in $A$.
The condition $\operatorname{reach}(A)>0$, introduced by Federer~\cite{fed59} is a
cornerstone in the geometric measure theory. This condition is a
considerable generalization of the notion of convexity [as the
convexity of $A$ is equivalent to $\operatorname{reach}(A)=\infty$ but a set
with $\operatorname{reach}(A)>0$ can be highly nonconvex].

\subsection{A Glivenko--Cantelli result for cone-convex sets}\label{GCresult}

We will next establish a GC-result for the class ${\mathcal C}_{\rho
,h}$ of nonempty compact sets fulfilling the $\rho,h$-cone-convex
condition (\ref{pzz}). In fact, we will establish the result for the
larger class ${\mathcal C}^U_{\rho,h}$ of \textit{closed} $\rho
,h$-cone-convex sets (thus we may drop the boundedness condition).
Since any closed convex set is in ${\mathcal C}_{\rho,h}^{U}$ we thus have
an extension of the well-known GC-type result for the class of closed
convex sets
(see, e.g., Talagrand~\cite{tal87}). Also, the following result provides a
strict generalization of Theorem~5 in Cuevas, Fraiman and Pateiro-L{\'
o}pez \cite{cue12} about the
GC-property for sets with reach $\geq h$: indeed, note that, as
established in that paper (Propositions 1 and 2), if a set fulfils
reach $\geq h$ then the outer $h$-rolling property holds and hence the
set belongs to the class ${\mathcal C}_{\rho,h}$ for $\rho<\pi/2$.
Moreover, the boundedness assumption is dropped here. Similar comments
and conclusions also hold for the class $\tilde{\mathcal C}^U_{\rho,h}$
of closed sets in ${\mathbb R}^d$ fulfilling the condition (\ref{rhoccc})
of cone-convexity by complements. All this is summarized in the
following result which will have some usefulness in the next section.

\begin{theorem}\label{GCtheorem} Let $\rho\in(0,\pi)$ and $h>0$. Let
$P$ be a probability measure absolutely continuous with respect to the
Lebesgue measure $\mu$. Then
\begin{longlist}[(a)]
\item[(a)]The class ${\mathcal C}^U_{\rho,h}$ of nonempty closed sets
fulfilling the $\rho,h$-cone-convex
condition (\ref{pzz}) is a $P$-uniformity class (and in particular a
$P$-Glivenko--Cantelli-class).
\item[(b)] The same conclusion holds for the class $\tilde{\mathcal
C}^U_{\rho,h}$ of closed sets fulfilling the condition (\ref{rhoccc}) of
cone-convexity by complements.
\end{longlist}
\end{theorem}

\begin{pf} (a) Let us first establish the result for the subclass $\A
\subset{\mathcal C}_{\rho,h}$ of sets in ${\mathcal C}_{\rho,h}$
included in a common compact set $K$.
From Proposition~\ref{muS}, both ${\mathcal C}_{\rho,h}$ and $\A$ are
$P$-continuous families. Given a sequence $\{A_n\}\subset\A$ there
exist (since $\A$ is relatively compact in the space of compact sets
endowed with the Hausdorff metric) a subsequence $\{A_{n_k}\}$
and a compact nonempty set $A$ such that $d_H(A_{n_k},A)\rightarrow
0$. From Theorem~\ref{teo2},
$A\in\C_{\rho,h}$ and hence $A\in\A$. From Theorem~\ref{teo1},
$d_H(\partial A_{n_k},\partial A)\rightarrow0$. Therefore, the class
$\partial\A=\{\partial A\dvtx A\in\A\}$ is compact in $\M$, thus
fulfilling the above mentioned sufficient condition in Billingsley and Tops{\o}e \cite{bil67}, Theorem~4.
This entails the $P$-uniformity property for the class $\A$.

Finally, given $\varepsilon>0$ take a large enough $R$ such that $P
(B(0,R)^c )< \varepsilon/8$. Let $K=\overline{B(0,R)}$.
If the weak convergence $P_n\stackrel{w}\longrightarrow P$ holds we
have, for large enough $n$, $P_n(K^c)<\varepsilon/4$. Then, denoting
${\mathcal C}^U_{\rho,h}={\mathcal D}$,
\begin{eqnarray*}
\sup_{A\in{\mathcal D}} \bigl|P_n(A)-P(A)\bigr |&\leq& \sup
_{A\in
{\mathcal D}} \bigl|P_n(A\cap K)-P(A\cap K) \bigr|\\
&&{}+\sup
_{A\in{\mathcal
D}} \bigl|P_n \bigl(A\cap K^c \bigr)-P
\bigl(A\cap K^c \bigr)\bigr |
\\
&\leq& \sup_{A\in{\mathcal D}} \bigl|P_n(A\cap K)-P(A\cap K) \bigr| +
P_n \bigl(K^c \bigr)+P \bigl(K^c \bigr)<
\varepsilon
\end{eqnarray*}
for $n$ large enough, since $A\cap K$ belongs to the class
$\A$.

(b) The result follows directly from (a) and from Proposition~\ref
{inclusiones}(d), which establishes that $\tilde{\mathcal C}_{\rho
,h}\subset{\mathcal C}_{\rho',h'}$ for suitable values of $\rho'$,
$h'$. Note also that, from Propositions \ref{inclusiones}(d) and \ref
{muS}, $\tilde{\mathcal C}_{\rho,h}$ is also a $P$-continuity class.
\end{pf}

\section{Estimation of cone-convex sets}\label{estimacion}

This section is devoted to the study of the asymptotic properties of
the two notions of cone-convex hull (when applied to a sample $\aleph
_n$) given in Definition~\ref{envolvente cono convexa}.

First, we obtain consistency and convergence rates for the $\rho,h$-cc
estimator ${\mathbb C}_{\rho,h}(\aleph_n)$. Second, we give convergence
rates for the $\rho,h$-ccc convex hull $\tilde{\mathbb C}_{\rho
,h}(\aleph
_n)$. Some key elements in the proof of the $\rho,h$-ccc case are the
notion of \textit{unavoidable families} (as in Pateiro-L{\'o}pez and Rodr{\'{\i
}}guez-Casal \cite{pat08}) and some
results on volume functions in Stach{\'o} \cite{sta76}.

\subsection{Consistency and rates for the cone-convex hull}\label{cons-rhocc}

The following consistency result is a direct consequence of our
GC-result (Theorem~\ref{GCtheorem}).

\begin{theorem}\label{GC-consistencia} Let $P$ be a
probability measure on ${\mathbb R}^d$, absolutely continuous with
respect to the Lebesgue measure $\mu$. Assume that $P$ has a compact
support $S$. Let $X_1,\ldots,X_n$ be a sample drawn from $P$. Denote
$\aleph_n=\{X_1,\ldots,X_n\}$.
\begin{longlist}[(a)]
\item[(a)] If $S$ is $\rho,h$-cone convex, then the sequence ${\mathbb
C}_{\rho
,h}(\aleph_n)$ of $\rho,h$-cone-convex hulls of $\aleph_n$ fulfills
%
\begin{equation}
d_H \bigl({\mathbb C}_{\rho,h}(\aleph_n),S \bigr)
\to0\qquad \mbox{a.s.}\quad \mbox{and} \quad d_\nu \bigl({\mathbb C}_{\rho,h}(
\aleph_n),S \bigr)\to0\qquad \mbox{a.s.}
\end{equation}
for any measure $\nu$, finite on compact sets, whose restriction to $S$
is absolutely continuous with respect to $P$.

\item[(b)] A similar result holds for the sequence $\tilde{\mathbb C}_{\rho
,h}(\aleph_n)$ of $\rho,h$-cone-convex hulls by complement, if we
assume that $S\in\tilde\C_{\rho,h}$.
\end{longlist}
\end{theorem}
\begin{pf}
(a) The first result, $d_H({\mathbb C}_{\rho,h}(\aleph_n),S)\to0$,
a.s. is obvious since $d_H(\aleph_n,S)\to0$ a.s. and
$\aleph_n\subset{\mathbb C}_{\rho,h}(\aleph_n)\subset S$.

As for the second result, note that
$d_\nu({\mathbb C}_{\rho,h}(\aleph_n),S)=\nu({\mathbb C}_{\rho
,h}(\aleph
_n)\setminus S)+\nu(S\setminus{\mathbb C}_{\rho,h}(\aleph_n))$.
The first term in the right-hand side is 0 a.s. As for the second one,
since $\nu$ is absolutely continuous with respect to $P$ on $S$ and $S$
is the support of~$P$, we only need to prove (from the well-known
$\varepsilon$--$\delta$ characterization of absolute continuity, when
$\nu$
is finite) that $P(S\setminus{\mathbb C}_{\rho,h}(\aleph_n))\to0$,
a.s. Indeed,
%
\begin{eqnarray}
P \bigl(S\setminus{\mathbb C}_{\rho,h}(\aleph_n)
\bigr)&=&P(S)-P \bigl({\mathbb C}_{\rho
,h}(\aleph_n) \bigr)
\\
\label{GCtasas}&\leq&\bigl|P(S)-{\mathbb P}_n \bigl({\mathbb C}_{\rho,h}(
\aleph_n) \bigr)\bigr|
\nonumber
\\[-8pt]
\\[-8pt]
\nonumber
&&{}+\bigl|{\mathbb P}_n \bigl({\mathbb
C}_{\rho,h}(\aleph_n) \bigr)-P \bigl({\mathbb
C}_{\rho,h}( \aleph_n) \bigr)\bigr|.
\end{eqnarray}
The first term is identically 0 a.s. The second one converges to 0 a.s.
from Theorem~\ref{GCtheorem}(a).

(b) The proof of (b) is completely analogous using Theorem~\ref
{GCtheorem}(b).
\end{pf}

\begin{remark}
A similar $d_\nu$-consistency result can be obtained by combining
Theorem~\ref{teo1} above
with Theorem~2 in Cuevas, Fraiman and Pateiro-L{\'
o}pez \cite{cue12}. However, Theorem~\ref{GC-consistencia}
provides a more
direct proof with an additional advantage: let us assume that $S$
belongs to a suitable subclass ${\mathcal A}\subset{\mathcal C}_{\rho
,h}$, and the estimator $S_n$ is chosen in that class. If the
convergence rate of $\sup_{A\in{\mathcal A}}|{\mathbb P}_n(A)-P(A)|$ is
known, then from (\ref{GCtasas}), the same convergence rate would
immediately apply to~$S_n$.
\end{remark}

The following theorem provides convergence rates in the Hausdorff
metric. Let us first recall (e.g., Cuevas and Fraiman \cite{cf97}) that (taking the
Lebesgue measure~$\mu$ as a reference) a set $S\subset{\mathbb R}^d$ is
said to be \textit{standard} with respect to a Borel measure~$\nu$ if
there exist $\lambda>0$, $\delta>0$ such that
%
\begin{equation}
\label{estandar} \nu \bigl(B(x,\varepsilon)\cap S \bigr)\geq\delta\mu \bigl(B(x,
\varepsilon) \bigr) \qquad\mbox{for all } x\in S, 0<\varepsilon\leq\lambda.
\end{equation}

\begin{theorem} \label{dHrates} Assume that $X_1,X_2,\ldots,X_n,\ldots
$ are
i.i.d. observations drawn from a distribution $P_X$ with support $S$.
Assume also that $S$ is compact and standard with respect to $P_X$.
Denote $\aleph_n= \{X_1,X_2,\ldots,X_n \}$. Then
\begin{longlist}[(a)]
\item[(a)] if $S\in\C_{\rho,h}$ then
$
d_H (\mathbb{C}_{\rho,h}(\aleph_n),S )= \mathcal{O} (
(\frac{\log n}{n} )^{1/d})$ a.s.

\item[(b)] The same conclusion holds for the estimator $\tilde{\mathbb C
}_{\rho,h}(\aleph_n)$ whenever the assumption $S\in\C_{\rho,h}$ is
replaced with $S\in\tilde\C_{\rho,h}$.
\end{longlist}
\end{theorem}
\begin{pf}(a)
Let us first consider the case $S\in\C_{\rho,h}$. Since $\aleph
_n\subset\mathbb{C}_{\rho,h}(\aleph_n)\subset\mathbb{C}_{\rho
,h}(S)=S$, the result follows directly from the following theorem given
in Cuevas and Rodr{\'{\i
}}guez-Casal \cite{cue04}, Theorem~3.

\begin{theorem*}
Let $X_1,X_2,\ldots$ be a sequence of i.i.d. observations drawn from a
distribution $P_X$ on ${\mathbb R}^d$. Assume that the support $S$ of
$P_X$ is compact and standard with respect to $P_X$. Then
\[
\limsup_{n\rightarrow\infty} \biggl(\frac{n}{\log n} \biggr)^{{1}/{d}}d_H(
\aleph_n,S)\leq \biggl(\frac{2}{\delta\omega_d} \biggr)^{{1}/{d}}\qquad
\mbox{a.s.},
\]
where $\omega_d$ is the Lebesgue measure of the unit ball in ${\mathbb
R}^d$ and $\delta$ is the standardness constant in (\ref{estandar}) for
$\nu=P_X$.
\end{theorem*}

(b) The proof for the case $S\in\tilde\C_{\rho,h}$ is identical.
\end{pf}

We will now study the rates of convergence for
$d_\mu(\mathbb{C}_{\rho,h}(\aleph_n),S )$,
with \mbox{$S\in\C_{\rho,h}$}.

We will need an assumption established in terms of the so-called
$t$-\textit{inner parallel set} of $S$, defined as
$S\ominus tB(0,1)= \{x\in S\dvtx B(x,t)\subset S\}$. The
inner parallel set appears as the result of applying the \textit{erosion
operator} $\ominus$ defined in the mathematical theory of image
analysis; see Serra \cite{ser82}. Also, the inner parallel set has received
some attention in differential geometry, on account of the regularity
properties of its boundary; see Fu \cite{fu85} and Remark~\ref{ipb} below.

\begin{theorem} \label{dmurates} Let $S\subset{\mathbb R}^d$
fulfilling the assumptions of Theorem~\ref{dHrates}. Moreover, let us
assume that
%
\begin{equation}
\label{muint} P_X \bigl(S\setminus S \ominus tB(0,1) \bigr)=
\mathcal{O}(t).
\end{equation}
Then,
\begin{longlist}[(a)]
\item[(a)] if $S\in\C_{\rho,h}$,
$d_{P_X} (\mathbb{C}_{\rho,h}(\aleph_n),S )=\mathcal{O}
( (\frac{\log n}{n} )^{1/d} )$ a.s.

\item[(b)] The same conclusion holds for the estimator $\tilde{\mathbb
C}_{\rho
,h}(\aleph_n)$ if we assume $S\in\tilde\C_{\rho,h}$.
\end{longlist}
\end{theorem}
\begin{pf} (a) Since
$d_{P_X} (\mathbb{C}_{\rho,h}(\aleph_n),S )=P_X (S\setminus
\mathbb{C}_{\rho,h}(\aleph_n) )$,
it suffices to show that
if $\eps_n=d_H(\aleph_n,S)$, then there exists $k\in{\mathbb R}$ such
that, with probability one, for~$n$ large enough,
%
\begin{equation}
\label{ec1} S\ominus k\eps_n B(0,1)\subset\mathbb{C}_{\rho,h}(
\aleph_n).
\end{equation}
Indeed, in this case we would have (using again Cuevas and Rodr{\'{\i
}}guez-Casal \cite{cue04}, Theorem~3)
\begin{eqnarray*}
P_X \bigl(S\setminus\mathbb{C}_{\rho,h}(\aleph_n)
\bigr)&\leq& P_X \bigl(S\setminus \bigl(S\ominus k
\varepsilon_n B(0,1) \bigr) \bigr)\\
&=&\mathcal{O}(k \varepsilon_n)=
\mathcal{O} \biggl( \biggl(\frac{\log n}{n} \biggr)^{1/d} \biggr)\qquad
\mbox{a.s.}
\end{eqnarray*}

More precisely, we will show that (\ref{ec1}) holds for
$k=3+2/\sin(\rho/2)$. Choose $n_0=n_0(\omega)$ such that
for $n>n_0$, $2\eps_n/\sin(\rho/2)<h/2$. Now, by contradiction if there
exists a sequence $x_n\in S\ominus k\varepsilon_nB(0,1)$ with $x_n\notin
{\mathbb C}_{\rho,h}(\aleph_n)$, we can find a sequence $B_n\in\C
_{\rho
,h}$ with $\aleph_n\subset B_n$ and $x_n\notin B_n$.

Since $\varepsilon_n=d_H(\aleph_n,S)$, there exists $X_i\in\overline
{B(x_n,\varepsilon_n)}$, but since $\aleph_n\subset B_n$ we also have
$\overline{B(x_n,\varepsilon_n)}\cap B_n\neq\varnothing$ and $\overline
{B(x_n,\varepsilon_n)}\cap B_n^c\neq\varnothing$. This entails the
existence of $z_n \in\partial B_n$, $z_n \in\overline{B(x_n, \eps
_n)}$. Since $B_n\in\C_{\rho,h}$, we may choose a unit vector $\xi_n$
with $C_{\rho,\xi_n,h}(z_n)\subset B_n^c$ which implies $C_{\rho,\xi
_n,h}(z_n)\cap\aleph_n=\varnothing$.
Let us now consider
$y_n=z_n+\frac{2\varepsilon_n}{\sin(\rho/2)}\xi_n$. If we prove $B(y_n,
2\varepsilon_n)\subset S$, we have got a contradiction with $d_H(S,\aleph
_n)=\eps_n$; indeed, from the definition of $y_n$ it is easy to see
that $B(y_n, 2\eps_n)\subset C_{\rho,\xi_n,h}(z_n)\subset B_n^c$ so, as
$\aleph_n\subset B_n$, one would have $B(y_n,2\varepsilon_n)\cap\aleph
_n=\varnothing$. Now, in order to prove $B(y_n, 2\varepsilon_n)\subset S$
recall that $B(x_n,k \varepsilon_n)\subset S$, so it suffices to check
$B(y_n,2\varepsilon_n)\subset B(x_n,k \varepsilon_n)$, but if $t\in
B(y_n,2\varepsilon_n)$, $\|t-x_n\|\leq\|t-y_n\|+\|y_n-z_n\|+\|z_n-x_n\|
\leq2\eps_n +2\varepsilon_n/\sin(\rho/2)+\varepsilon_n= \varepsilon
_n(3+2/\sin
(\rho/2))$.

(b) The result follows from the above conclusion (a) and Proposition~\ref{inclusiones}(d), since according to this result
$\tilde{\mathbb C}_{\rho,h}(\aleph_n)\in\C_{\rho',h'}$ so that
\[
{\mathbb C}_{\rho',h'}(\aleph_n)\subset\tilde{\mathbb
C}_{\rho
,h}(\aleph_n).
\]
\upqed\end{pf}
%

\begin{remark}\label{tasas-conv} The convergence order obtained in
Theorem~\ref{dmurates} is the same found in D{\"u}mbgen and Walther \cite{dum96} for the case
in which the ordinary notion of convexity for $S$ (and the
convex hull for the estimator) are used, instead of the
much more general concept of cone-convexity considered here.
The same behavior is found in Rodr{\'{\i}}guez-Casal \cite{rod07} for the intermediate case
in which $r$-convexity is assumed.
\end{remark}

\begin{remark}\label{ipb} Note that
$S\setminus S \ominus tB(0,1)$ is the set of points in $S$ within a
distance from $\partial S$ smaller than $t$. Thus, condition (\ref
{muint}) has a clear intuitive interpretation, connected with some key
concepts in Geometric Measure Theory.
To begin with, let us recall that the erosion operator $\ominus$
provides (as well as the dual dilation operator $\oplus$) a well-known
standard ``smoothing'' procedure in the mathematical theory of image analysis.
Now, to give a more precise interpretation of condition (\ref{muint})
let us assume that $P_X$ is uniform, that is, proportional to the
Lebesgue measure $\mu$ (similar conclusions can be drawn when $P_X$
fulfils $c_1\mu\leq P_X\leq c_2 \mu$ for some constants $c_1,c_2>0$).
Note that, if we denote $T=S^c$, we have $S\setminus S \ominus
tB(0,1)\subset B(T,t)\setminus T$.
We thus have that (\ref{muint}) will hold whenever $\mu
(B(T,t)\setminus
T)=\mathcal{O}(t)$.
A sufficient condition for this would be the celebrated Federer's
\textit{positive reach condition}, a geometric smoothness notion introduced
at the end of Section~\ref{GCsum} above.
More specifically, it is proved in Federer \cite{fed59}, Theorem~5.6, that if
$\operatorname{reach}(\bar T)=R$, then $\mu(B(T,t)\setminus T)$ is a polynomial
in $t$, of degree $d$, for $t\in[0,R)$; in particular, (\ref{muint})
holds. Also, the finiteness of the outer Minkowski content of $T$
(defined by $L_1=\lim_{t\to0}\mu(B(T,t)\setminus T )/t$; see
Ambrosio, Colesanti and Villa \cite{amb08}) is a sufficient condition for (\ref{muint}).
\end{remark}

The following result shows that the boundary of $S$ can be estimated as
well, with rates of the same order, under our cone-convexity assumption.

\begin{corollary} \label{bordes} Under the assumptions of Theorem~\ref
{dmurates}\textup{(a)}, we have that, with probability one, for $n$ large enough,
$d_H (\partial S,\partial\mathbb{C}_{\rho,h}(\aleph_n) )\leq
kd_H(S,\aleph_n)$,
where
$k= (3+2/\sin(\rho/2) )$. A similar result holds for the $\rho
,h$-ccc-hull $\tilde{\mathbb C}_{\rho,h}(\aleph_n)$ if we assume the
conditions of Theorem~\ref{dmurates}\textup{(b)}. In this case, the constant
$k$ is replaced with $k'= (3+2/\sin(\rho'/2) )$, where $\rho'$ is
the angle defined in Proposition~\ref{inclusiones}\textup{(d)}.
\end{corollary}
\begin{pf} In the case (a), the result follows from the content
relation (\ref{ec1}) together with $\partial{\mathbb C}_{\rho
,h}(\aleph
_n)\subset S\setminus\operatorname{int} (S\ominus k\varepsilon_nB(0,1)
)\subset B(\partial S,k\varepsilon_n)$ and the fact that, for any $x\in
\partial S$
there is a sample point $X_i$ such that
$\Vert x-X_i\Vert\leq\varepsilon_n$. Then,
in the segment joining $x$ and $X_i$ there must be necessarily a point of
$\partial{\mathbb C}_{\rho,h}(\aleph_n)$.
In the case (b), we use again $\mathbb{C}_{\rho',h'}(\aleph
_n)\subset
\tilde{\mathbb C}_{\rho,h}(\aleph_n)$ and Proposition~\ref
{inclusiones}(d).
\end{pf}
%

\subsection{Convergence rates in mean. Unavoidable families}\label
{rates-rhoccc}

We now focus on the convergence rates
for the ``mean error in measure'' ${\mathbb E}d_{P_X}(\tilde{\mathbb
{C}}_{\rho,h}(\aleph_n),S)$
where $P_X$ denotes the distribution of $X$. As we will see, the corresponding
proof will involve some interesting methodological differences with the
techniques used so far. In particular, relying on some ideas in
Pateiro-L{\'o}pez and Rodr{\'{\i
}}guez-Casal \cite
{pat12}, we will use the auxiliary notion of \textit{unavoidable families
of sets} which is next introduced and analyzed.
Under suitable conditions ensuring $\C_{\rho,h}\subset\tilde\C
_{\rho',h'}$,
it should be possible also to obtain an analogous result for the cc-hull
estimator ${\mathbb{C}}_{\rho,h}(\aleph_n)$. However, this technical
issue will not be considered here.

\textit{Unavoidable families of sets}.
Given $a\in(0,\pi)$ and $b>0$, denote
${\mathcal G}_{a,b}$ the family of all cones with opening
angle $a$ and height $b$, that is, ${\mathcal G}_{a,b}= \{C_{a,\xi
,b}(x)\dvtx\break  x\in{\mathbb R}^{d}, \Vert\xi\Vert= 1 \}$.

\begin{definition}\label{inevitables} A family of nonempty sets
${\mathfrak U}$ is said to be \textit{unavoidable} for another family
of sets $\Lambda$ if for
each $L\in\Lambda$ there exists $U\in{\mathfrak U}$ with $U\subset L$.
\end{definition}

The reason for using this notion here is as follows. Let
$\Lambda_{\rho,h}(x)=\{C\in{\mathcal G}_{\rho,h}\dvtx x\in C\}$, that is,
$\Lambda_{\rho,h}(x)$ is the family of $\rho,h$-cones which
include the point $x$. Assume that we are able to find for each $x
\in S$ a suitable finite family ${\mathfrak U}_{x,\rho,h}$,
unavoidable for $\Lambda_{\rho,h}(x)$. Assume
also that $X$ has a density $f$ satisfying $0<k_1 \leq f(x) \leq
k_2< \infty$ for almost all $x$ in $S$. We would then have
\begin{eqnarray}\label{equnav}
\mathbb{P} \bigl(x\in S\setminus\tilde{\mathbb{C}}_{\rho,h} (
\aleph_n) \bigr)&= &\mathbb{P} \bigl(\exists C\in\Lambda_{\rho,h}(x)
\dvtx C\cap\aleph_n= \varnothing \bigr) \nonumber\\
&\leq&\sum
_{U\in\mathfrak{U}_{x,\rho,h}} \mathbb{P}(U\cap\aleph_n=\varnothing )\quad
\mbox{and}
\nonumber\\
\mathbb{E} \bigl(d_{P_X} \bigl(S,\tilde{
\mathbb{C}}_{\rho,h}(\aleph_n) \bigr) \bigr) &=& \mathbb{E}\int
_S \ind_{\{x\in
S\setminus\tilde{\mathbb{C}}_{\rho,h} (\aleph_n)\}}f(x)\,dx
\nonumber
\\[-8pt]
\\[-8pt]
\nonumber
&=&\int_S
\mathbb{P} \bigl(x\in S\setminus\tilde{\mathbb{C}}_{\rho,h} (
\aleph_n) \bigr)f(x)\,dx\\
&\leq& k_2 \int_S \sum
_{U\in\mathfrak{U}_{x,\rho,h}} \mathbb{P}(U\cap\aleph_n=\varnothing )\,dx
\nonumber\\
&\leq& k_2 \int_S \sum
_{U\in\mathfrak{U}_{x,\rho,h}} \bigl(1-k_1 \mu(U\cap S)
\bigr)^n \,dx,\nonumber
\end{eqnarray}
where in the last inequality we have also used
that $f$ is bounded from below.
So, in order to find rates of convergence for ${\mathbb E}
(d_{P_X} (S, \tilde{\mathbb{C}}_{\rho,h}(\aleph_n) ) )$ the
problem can be reduced to find, for each $x \in S$, a finite
unavoidable family ${\mathfrak U}_{x,\rho,h}$ such that $k_1 \mu(U
\cap
S)$ is large enough. Such families are described in the following
proposition whose proof is given in the \hyperref[app]{Appendix}.

\begin{proposition} \label{propunav}
Let $\gamma= \rho$ if $\rho\leq\pi/3$ and $\gamma= (\pi- \rho)/2$
otherwise.
Take $h_1 = \frac{h}{2} \sin(\frac{\rho}{2})$. Given $x \in S$,
consider a
minimal covering of the closed ball $ \overline{B(x,h_1)}$ with closed
cones of angle $\gamma/2$, axis $\xi_j$ and height
$h_1$, $\{ \overline{C_{\gamma/2, \xi_j, h_1}(x)},\break   \Vert\xi_j
\Vert
=1,  j=1, \ldots, k \}$. Then the family
\[
{\mathfrak U}_{x, \rho, h} = \bigl\{C_{\gamma/2, \xi_j, h_1}(x), \Vert
\xi_j \Vert=1, j=1, \ldots k \bigr\},
\]
is unavoidable for $\Lambda_{\rho, h}(x)$. Moreover, the cardinality of
${\mathfrak U}_{x, \rho, h}$ does
not depend on $x$.
\end{proposition}

We now establish the main result of this section. Again the proof is
given in the \hyperref[app]{Appendix}.

\begin{theorem}\label{thunav} Let $S\subset{\mathbb R}^{d}$, $S\in
\tilde{\mathcal C}_{\rho,h}$ and, for $z\geq 0$, $F(z)=\mu(\{x\in S\dvtx\break   d(x,\partial
S)\leq z\} )$. Assume that $F'$ is bounded in a neighborhood of 0, and
$X_1,X_2,\ldots$ are i.i.d. drawn from a
distribution $P_X$ with support $S$. Let us suppose that $P_X$ is
absolutely continuous with $\mu$-density $f$ such that $0 < k_1 \leq
f(x) \leq k_2 < \infty$ for some constants $k_1,k_2$ and for
almost all $x\in S$. Then
$\mathbb{E} (d_{P_X} (S,\tilde{\mathbb{C}}_{\rho,h}(\aleph_n)
) )=\mathcal{O}(n^{-1/d})$.
\end{theorem}

\begin{remark} \label{rrra}
(a) Note that for any given compact set $S\subset{\mathbb R}^d$
with $\mu(\partial S)=0$, $F'(0)$ the outer Minkowski content of $S^c$ (defined by $\lim_{\varepsilon\to0}\frac{\mu(B(S^c,\varepsilon)\setminus S^c)}{\varepsilon
}$). See Ambrosio, Colesanti and Villa \cite{amb08} for a deep study on this notion.

(b) The rate of convergence we have obtained is slower (when
$d=2$) than the one obtained in Pateiro-L{\'o}pez and Rodr{\'{\i
}}guez-Casal \cite{pat12} (Theorem~1) for
$r$-convex sets in ${\mathbb R}^2$, fulfilling a double rolling
condition. In return, the class of cone-convex sets we are considering
is much larger and we have no restriction on the dimension.
\end{remark}

\section{A stochastic algorithm for ccc-hulls}\label{algoritmo}

We offer here a relatively simple stochastic algorithm to approximately
calculate the cone-convex hull by complement, $\tilde{\mathbb C}_{\rho
,h}(\aleph_n)$ for a given random sample $\aleph_n=\{X_1,\ldots,X_n\}$.

As explained in Section~\ref{aplicaciones}, $\tilde{\mathbb
C}_{\rho
,h}(\aleph_n)$ is a close analogue of the $r$-convex hull previously
considered in the literature,
%
\begin{equation}
r\operatorname{conv}(\aleph_n)=\bigcap_{\{y: B(y,r)\cap\aleph_n=\varnothing
\}
}B(y,r)^c.\label{cierre-conv}
\end{equation}
An exact algorithm for the calculation of (\ref{cierre-conv}) for
samples in ${\mathbb R}^2$ can be found in the R-package
\textit{alphahull}, described in Pateiro-L{\'o}pez and Rodr{\'{\i
}}guez-Casal~\cite{pat10}.

The numerical treatment of the ccc-hull $\tilde{\mathbb C}_{\rho
,h}(\aleph_n)$ is a bit harder. This is essentially due to the lack of
rotational symmetry of the ``primary blocks'' used in the construction
of $\tilde{\mathbb C}_{\rho,h}(\aleph_n)$, which are finite cones,
instead of the balls of (\ref{cierre-conv}).

Our algorithm is based on the insightful heuristic description of (\ref
{cierre-conv}) given in Edelsbrunner and M{\"u}cke \cite{edel94}:
\textit{``Think of ${\mathbb R}^3$ filled with styrofoam and the points in
$\aleph_n$ made of more solid material such as rock. Now imagine a
spherical eraser with radius $r$. It is omnipresent in the sense that
it carves out styrofoam at all positions where it does not enclose any
of the sprinkled rocks, that is, points of $\aleph_n$. The resulting
object will be called the $r$-hull.''}

In our case, the ``eraser element'' is a finite cone
$C_{\rho,h,\xi}(x)$ instead of a ball $B(x,r)$. So, in order to move
the eraser we should in fact vary two parameters: the vertex $x$ and
the axis direction $\xi$ (since the angle $\rho$ and the height $h$
remain fixed).

Our proposal is essentially based on the idea of choosing these two
parameters with an ``oriented random procedure'': we pick up randomly
the vertex $x$ and then we erase as much styrofoam as possible by
rotating the cone $C_{\rho,h,\xi}(x)$ for all directions $\xi$ with
$C_{\rho,h,\xi}(x)\cap\aleph_n=\varnothing$. For $\theta\in[0,\pi/2]$,
let us denote $R^x_\theta(u)$ the clockwise rotation of angle $\theta$
with center in $x$ of the vector $u$, (if $\theta\in[-\pi/2,0)$ we
take the counter clockwise rotation). Then our algorithm is, in $\R
{2}$, as follows:

\noindent\hrulefill
\begin{itemize}
\item[] INPUT: A sample
$\aleph_n=\{X_1,\ldots,X_n\}\subset{\mathbb R}^2$, the cone parameters
$\rho\in(0,\pi]$ and $h>0$, a rectangle $E=[a,b]\times[c,d]$ with
$\aleph_n\subset E$, $N$ a large positive integer indicating the number
of full iterations of steps 1--3 below.
\item[] STEP 1. \emph{Generating random cones}: Choose at random a cone
vertex $x\in E$ and a cone axis $\xi$ with $\Vert\xi\Vert=1$ and
consider the cone $C_{\rho,h,\xi}(x)$.
\item[] STEP 2. \emph{Checking for an empty cone}: If
$C_{\rho,h,\xi}(x)\cap\aleph_n\neq\varnothing$ go back to step 1.
\item[] STEP 3. \emph{Erasing a maximal cone}: If
$C_{\rho,h,\xi}(x)\cap\aleph_n=\varnothing$ erase the maximal cone with
vertex $x$ not containing any sample point. That is, find
\begin{eqnarray*}
\theta_0& =& \max_{\theta\in[0,\pi/2]} \bigl\{\theta\dvtx
C_{\theta
/2,h,R^x_{\theta/2}(\xi)}(x)\cap\aleph_n=\varnothing \bigr\},
\\
\theta_1 &=& \min_{\theta\in[-\pi/2,0)} \bigl\{\theta\dvtx
C_{\theta
/2,h,R^x_{\theta/2}(\xi)}(x)\cap\aleph_n=\varnothing \bigr\}.
\end{eqnarray*}
Then erase the $h$-cone $C$ with vertex $x$ and sides of length $h$
along the directions $R^x_{\theta_1}(\xi)$ and $R^x_{\theta_0}(\xi)$.
That is, replace $E$ with $E\setminus C$ and go back to step 1.
\item[] OUTPUT: The set $E$ resulting after step 3 has been performed
$N$ times. So, $N$ is the number of erasing cones
during the iteration process.
\end{itemize}
\hrulefill\vspace*{6pt}

\textit{Some comments on the algorithm}.
\begin{longlist}[1.]
\item[1.] The \texttt{R} code of this algorithm (including detailed comments)
can be~down\-loaded from \url
{http://www.uam.es/antonio.cuevas/exp/ccc-algorithm.txt}.

\item[2.] The accuracy of the algorithm could be improved with some simple
changes. For example, we might choose the vertices in step 1 with a
probability measure whose density is inversely proportional to a kernel
density estimator of the underlying distribution of the sample. Of
course, the idea is to increase the probability of selecting vertices
in ``empty areas.'' We might also improve the efficiency by using the
convex hull (or the $h$-convex hull) of the sample as the initial
``frame'' $E$ to draw the cones. However, we have omitted such
modifications in order to present the idea in the most simplest way.

\item[3.] Finding exact (nonstochastic) algorithms to calculate both $\tilde
{\mathbb C}_{\rho,h}(\aleph_n)$ and ${\mathbb C}_{\rho,h}(\aleph
_n)$ is
a much harder problem, far beyond the scope of this paper. The exact
calculation of ${\mathbb C}_{\rho,h}(\aleph_n)$ seems particularly
difficult. The trouble lies in the fact that the cc-property, similarly
to the analogous ``outer sphere'' or ``rolling-ball property,'' does not
seem to provide a ``canonical way'' to construct a small enough set
including the sample points and fulfilling the cc-property. On the
contrary, the definition of the ccc-property implicitly includes a
mechanism to construct the ccc-hull.

\item[4.] We present here the algorithm for the two-dimensional case
$d=2$ since this is, by far, the most important case in the usual
applications and the presentation becomes a bit simpler. However,
the algorithm can be extended, with no essential change, to $d=3$
and, in fact, the basic idea would also work for $d>3$.

\item[5.] To give just an approximate idea of the execution time of our
algorithm, let us point out that the mean execution time over 1000 runs
(with $n=500$, $N=300$ cones $h=1/4$, $\rho=\pi/4$) was $36.453$
seconds for the set in the first example of Section~\ref{sec-examples} below. The corresponding standard desviation was $3.362$
seconds. We have used a processor Intel i7-2620M.
\end{longlist}

\begin{figure}[b]

\includegraphics{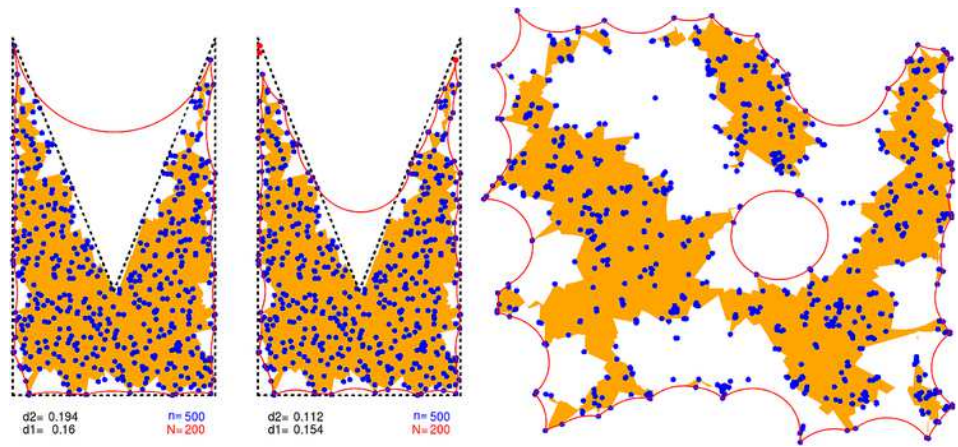}

\caption{Comparison of the ccc-hull (shaded area) and the $r$-convex
hull (boundary made of $r$-arcs).}
\label{fig4}
\end{figure}

\section{Some numerical results}\label{numerical-sec}

\subsection{Three examples}\label{sec-examples}
Just in order to gain some insight on the behavior of our
ccc-estimator we show here three examples. In all of them, we have
compared the ccc-hull with the above commented $r$-convex hull
(see, e.g., Pateiro-L{\'o}pez and Rodr{\'{\i
}}guez-Casal \cite{pat10} and references therein) which appears to
be the most direct competitor, as a generalization of the ordinary
convex hull.

The first example (left-hand and central panel of Figure~\ref{fig4}) illustrates the estimation (from $n=500$ uniform
points) of the $\pi/4$-cone convex set $S=[0,1]\times
[0,t+1/2]\setminus T$ where $t=\frac{1}{2}\tan(3\pi/8)$, $T$ being
the isosceles triangle with vertices $(0,t+1/2)$, $(1,t+1/2)$,
$(1/2,1/2)$. For the ccc-hull (the shaded area in the figures), we
have used $\rho=\pi/4$ and $h=1/2$ with $N=200$ cones. For the
$r$-convex hull (whose boundary is marked in continuous lines as a
union of $r$-circumference arcs), we took $r=1/2$ (left-hand
panel) and $r=1/4$ (central panel). The whole point of choosing
this set is to show that, even in very simple cases, the presence
of an inward nonsmooth peak can lead to a situation for which the
$r$-convex hull provides an ``oversmoothed'' estimation since the
estimator just cannot ``go inside'' the sharp ``gulf'' in the set.
This is not the case of the ccc-hull which is designed to deal
with such unsmooth situations. Of course, we might improve things
by choosing a smaller value of $r$ but, in any case, the
$r$-convex is inconsistent for any $r$ and, at the end, it will we
outperformed by the $\rho,h$-ccc hull, provided that a suitable
value of $\rho$ ($\leq\pi/4$ in this case) is chosen.

\begin{figure}

\includegraphics{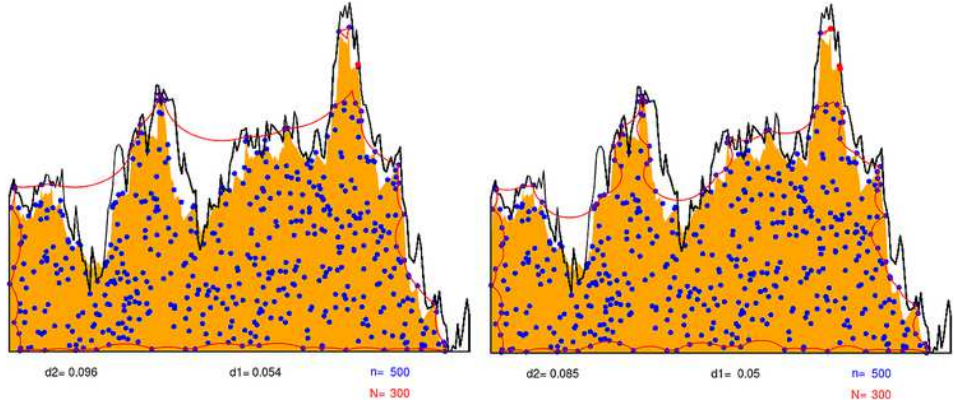}

\caption{Hypograph of a trajectory of a Brownian motion: the
shaded area is the ccc-hull;
continuous lines made of $r$-arcs correspond to the boundaries
of the $r$-convex hull for $r=1/8$ (left) and $r=1/16$ (right).}
\label{figura6}
\end{figure}

The second example (right-hand panel of Figure~\ref{fig4}) shows
the behavior of our estimator (with $\rho=\pi/3$, $h=1/8$) when
compared with the $r$-convex hull (the circumference arcs in
continuous lines) with $r=1/10$ for a sample of points that
represent the locations $(x,y)$ of bramble canes in a field of 9
square meters, rescaled to the unit square. This data set can be
found in the R-library \textit{spatstat};
see \cite{hut79} and \cite{dig83} for further details;
we have ignored the labels
identifying different classes of plants, according to their ages.
Of course, in this example there is no ``true'' set to give an
objective comparison. We can see, however, how both estimators give
a quite different estimation of the ``habitat'' of these plants
and the ccc-hull seems better adapted to detect the absence of
canes in some areas.
Finally, the third example shows the estimation of a quite
irregular set: the hypograph of the trajectory of a Brownian
motion on the unit interval. We define the hypograph of a positive
function $f$ defined on $[a,b]$ by $H(f)=\{(x,y)\dvtx x\in[a,b], 0\leq
y\leq f(x)\}$. In our example, the Brownian trajectory has been
shifted vertically in order to take all values above zero. The
estimation of hypographs is a major aim in the problem of
\textit{efficient boundary}, a relevant topic in econometrics; see, for
example, Simar and Wilson \cite{sim00}. An additional interest of this example is to show
how our ccc-hull can be adapted to incorporate the information
that our target set is an hypograph; this can be made by just
choosing vertical cones and restricting their rotation angle in
the algorithm. In this case, we have taken $n=500$, $\rho=\pi/6$,
$h=1$, with $N=300$ cones but the rotation angles in the algorithm
have been restricted between $5\pi/12$ and $7\pi/12$ in order the
keep the structure of an hypograph; see Figure~\ref{figura6}. The
parameters for the
$r$-convex estimator are $r= 1/8$ (left panel in Figure~\ref{figura6}),
and $1/16$ (right panel);
note that there is no way to adapt the $r$-convex hull to the
hypograph shape. In this case, the ccc-hull, with the hypograph
information incorporated, clearly outperforms the $r$-convex hull.

\subsection{Simulation outputs}
\begin{figure}[b]

\includegraphics{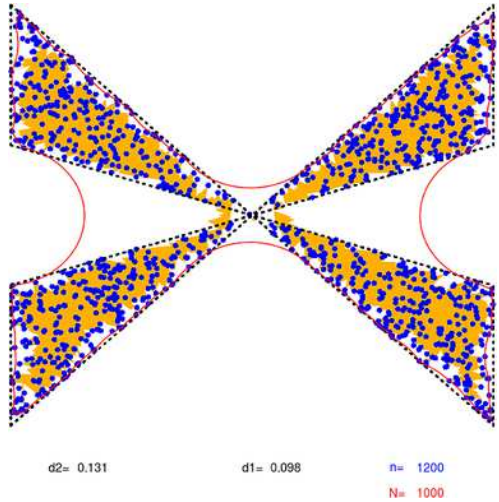}

\caption{Ccc-hull (shaded area) and $r$-convex hull (boundary
made of $r$-arcs) constructed from 1200 points of $S_1$. Here, $\rho
=\pi
/5$, $h=1/2$, $r=1/6$.} \label{Rplot}
\end{figure}

We have carried out a small simulation study to compare the
performance of the ccc-hull with that of the $r$-convex hull for
different sample sizes and\vadjust{\goodbreak} values of the parameters.
The target set $S_1=[0,1]^2\setminus
\bigcup_{i=1}^4T_i$
where the $T_i$ are triangles with vertices $(0,1)$, $(1/2,1/2)$,
$(1,1)$; $(0,0)$, $(1/2,1/2)$, $(1,0)$;
$(0,1/3)$, $(1/2,1/2)$, $(0,2/3)$ and $(1,1/3)$,\break  $(1/2,1/2)$, $(1,2/3)$.
This set is $\rho_0=2\arctan(1/3)$-cone convex. Figure~\ref{Rplot}
corresponds
to the case $\rho=\pi/5$, $h=1/2$, $N=1000$ and $r=1/6$ for the
$r$-convex hull with $n=1200$ uniform points.

Table~\ref{testM4} shows the expected
values, over 500 runs (and their standard deviations in
parenthesis) for the errors in measure [$d_1=d_\mu(S,{\mathbb
C}_{\rho,h}(\aleph_n) )$,
$d_2=d_\mu(S,C_{r}(\aleph_n) )$] of both estimators (the
ccc-hull and the $r$-convex hull), with different values
of the parameters $\rho$, $r $ and $h$. For
small sample sizes (200 in Table~\ref{testM4}) the $r$-convex hull
has a smaller error in measure. However, as the sample size
increases, (from 400 on the ccc-hull
outperforms the $r$-convex hull. We
have taken $N=200$ cones for the simulation, and the distances
were calculated by the Monte Carlo method using 4000 uniform
random observations.
%
\begin{table}
\caption{Average errors (and standard deviations) over 500 runs
for the ccc-hull ($d_1$) and the $r$-convex hull ($d_2$) in the
estimation of $S_1$} \label{testM4}
\begin{tabular*}{\textwidth}{@{\extracolsep{\fill}}lcccc@{}}
\hline
& \multicolumn{1}{c}{$\bolds{d_1}$} & \multicolumn{1}{c}{$\bolds{d_2}$} & \multicolumn{1}{c}{$\bolds{d_1}$} & \multicolumn{1}{c@{}}{$\bolds{d_2}$}\\[-6pt]
& \multicolumn{1}{c}{\hrulefill} & \multicolumn{1}{c}{\hrulefill} & \multicolumn{1}{c}{\hrulefill} & \multicolumn{1}{c@{}}{\hrulefill}\\
$\bolds{n}$ & $\bolds{\rho=\rho_0,h=1/3}$ & $\bolds{r=1/4}$ &$\bolds{\rho=\pi/5,h=1/2}$ & $\bolds{r=1/6}$\\
\hline
\phantom{0}200 & 0.204 (0.011) & 0.191 (0.009) & 0.197 (0.011) & 0.161 (0.010)\\
\phantom{0}400 & 0.138 (0.009) & 0.180 (0.008) & 0.134 (0.010) & 0.140 (0.008) \\
\phantom{0}600 & 0.107 (0.008) & 0.174 (0.007) & 0.105 (0.008) & 0.132 (0.007)\\
\phantom{0}800 & 0.090 (0.007) & 0.172 (0.007) & 0.089 (0.007) & 0.127 (0.007)\\
1000 & 0.080 (0.007) & 0.170 (0.007) & 0.078 (0.007) & 0.124 (0.006) \\
1200 & 0.070 (0.006) & 0.169 (0.006) & 0.070 (0.006) & 0.122 (0.006)\\
\hline
\end{tabular*}
\end{table}

%

\section{Final remarks: Some suggestions for further work}\label{final-sec}

In our view, the study of the following topics might be of interest in
connection with the notion of cone-convexity introduced in this paper.

\textit{Applications to home-range estimation}.
As commented above, our cone-convex hulls are in fact a considerable
generalization of the simpler classical notion of convex-hull. Such
generalizations (the $r$-convex hull is another example of them) are
relevant in those application fields where more flexible set estimators
are needed.
An example arises in zoology and ecology, in the problem of \textit{home
range estimation}. A commonly cited definition of animal's home range
is that of Burt~\cite{bur43}: ``that area traversed by the individual in
its normal activities of food gathering, mating and caring for young.''
The problem of estimating the home range from ``sightings'' or GPS
records of animal positions has received a considerable attention (see,
e.g., Anderson \cite{and82} for an introduction). As pointed out by Burgman and
Fox~\cite{bur03}, \textit{``Minimum convex polygons (convex hulls) are an
internationally accepted, standard method for estimating
species' ranges, particularly in circumstances in which presence-only
data are the only kind of spatially
explicit data available''}. These authors also discuss the obvious
drawbacks of the convex hull, and analyze in some detail the so called
$\alpha$-hulls (conceptually related with the $r$-convex hulls
discussed above) as a useful more flexible alternative. In fact, the
idea of considering different nonparametric estimators in home range
estimation is far from new. Many highly cited papers (Worton \cite{wor89},
Getz and Wilmers \cite{get04}, etc.) have considered this topic. Some of them, in
particular, Worton \cite{wor89}, analyze the use of auxiliary density
estimators to construct home range estimators. We believe that our
proposal here, based on the cone-convex hull, could be seen as a
further step in this advance toward flexibility and generality from the
classical approach based on the ``hull principle.'' The reason is that
our estimator could be suitable for those problems where highly
irregular shapes, including central holes of sharp inward peaks, are to
be expected, due to existence of geographic obstacles leading to
irregular habitats. For example, Getz and Wilmers \cite{get04} have suggested (in a
nonmathematical journal) an interesting class of estimators based on
the union of convex hulls of the nearest neighbours of every
sample point. These authors convincingly motivate their proposal on the
basis of detailed examples. Again, the point is the need of flexible,
general estimators for home range studies and related problems.
However, to the best of our knowledge, the theoretical properties of
that class of estimators have not been analyzed so far. In a way, our
proposal in this paper, aims at same goals having still in mind the
idea of extending the classical convex hull. While the detailed
analysis of such practical applications is beyond the scope of this
paper, we hope that the real-data example (not in zoology but in
botany) outlined in the previous section could give a hint on the
possible advantages of our estimators.

\textit{Inference on the parameter $\rho$}.
In our cone-convexity definitions, the parameter $\rho$ has an obvious
intuitive interpretation (in terms of the sharpest inward peak in the
domain), even more direct than that of the parameter $r$ in the
$r$-convexity property. So, given a domain $S$, the inference on the
largest value of $\rho$ fulfilling the cone-convexity property (for a
given $h$) might be of some interest from the image analysis point of
view. In particular, the study of a suitable test for the hypothesis
$H_0\dvtx \rho\geq\rho_0$ seems a natural aim. Note that in the case
$\rho
_0=\pi$ this would essentially amount to test convexity. The theory of
multivariate spacings, as developed, for example, by Janson \cite{jan87},
seems to be a relevant auxiliary tool
in this problem.

\textit{Cone-convexity for functions}.
Our cone-convexity concepts have been primarily defined for sets but
they could be extended in a natural way for real functions
$f\dvtx [a,b]^d\rightarrow{\mathbb R}$: we could say that $f$ is $\rho
$-cone-convex when the hypograph $H(f)=\{(x,y)\dvtx x\in[a,b]^d, y\leq
f(x)\}$ is $\rho$-cone-convex. The distance between two $\rho$-cc
functions might then be defined in terms of the Hausdorff distance
between the corresponding hypographs; similar ideas have been
considered elsewhere, for example, Sendov \cite{sen90}. On the one hand, this
Hausdorff-based metric would provide a ``visual'' proximity criterion
(potentially meaningful in many real-world applications) between the
data. On the other hand, the cone-convexity assumption would lead to a
natural way for data
smoothing. For example, in the setting of a nonparametric regression
model $y_i=f(x_i)+\varepsilon_i$ (with $d=1$), we could think of
recovering the function $f$ from the data $(x_i,y_i)$ under the
assumption that $f$ is $\rho$-cone-convex. Other applications to
Functional Data Analysis (in particular to supervised functional
classification) are also under study.

\begin{appendix}
\section*{Appendix}\label{app}

\begin{pf*}{Proof of Proposition~\ref{propunav}}
Let $C$ be a member of the class $\Lambda_{\rho, h}(x)$. Without loss
of generality, take
$C=C_{\rho,e_1,h}(0)$ where $e_1$ is the first vector in the canonical
basis. The reasoning for any other cone
$C=C_{\rho,\xi,h}(z)$ is reduced to this case by translation and/or rotation.

\begin{figure}

\includegraphics{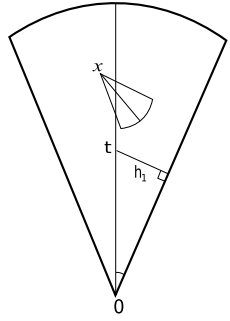}

\caption{Case 3, $C_{\gamma,\nu,h_1}(x) \subset C_{\rho,e_1,h}(0)$.}
\label{cono}\vspace*{-2pt}
\end{figure}

By definition of the class, we have $x\in C_{\rho,e_1,h}(0)$. As a
first step in the proof, it will be useful to consider three possible
situations regarding the position of $x$. In all three cases, we will
be able to find a cone \mbox{$C'=C_{\gamma,v,h_1}(x)\subset C$}.

\begin{longlist}[1.]
\item[1.] If $x\in C_{\rho,e_1,h/2}(0)$, then $C'\subset C_{\rho,e_1,h}(0)$
for $C'=C_{\gamma,e_1,h_1}(x)$.
\item[2.] If $x=re_1$ with $h/2\leq r < h$, then $C'\subset C_{\rho
,e_1,h}(0)$ for $C'=C_{\gamma,-e_1,h_1}(x)$ since $h_1$ is smaller
than the
distance from $he_1/2$, the middle point of the axis of $C=C_{\rho
,e_1,h}(0)$, to the boundary of $C$.
\item[3.] For any other $x \in C_{\rho,e_1,h}(0)$, note that $t=he_1/2$
corresponds to the least favorable position of the vertex of a $\gamma
,h_1$-cone, $C^\prime$ in order to get $C^\prime\subset C$. Let
$\lambda
_1$ be the solution of $ \langle\frac{\lambda_1t-x}{\|\lambda_1t-x\|
},\frac{-x}{\|x\|} \rangle= \cos(\gamma/2)$. Then if we take
a cone $C'=C_{\gamma,v,h_1}(x)$ with axis $v=\frac{\lambda_1 t-x}{\|
\lambda_1 t-x\|}$ we also get $C'\subset C_{\rho,e_1,h}(0)$; see Figure~\ref{cono}.
\end{longlist}

Finally, the unavoidable family is constructed by selecting a finite
number of directions $\xi_j$ such that the cones
$\{\overline{C_{\gamma/2,\xi_j,h_1}(x)}\}_{j=1,\ldots,k}$ are a minimal
covering of $\overline{B(x,h_1)}$. Indeed, given $C\in{\Lambda
}_{\rho
,h}(x)$ the point $x\in C$ is in one of the three previously considered
cases so that there exists a unit vector $\xi$ for which $C_{\gamma
,\xi
,h_1}(x)\subset C$.
Since $\{\overline{C_{\gamma/2,\xi_j,h_1}(x)}\}_{j=1,\ldots,k}$ is a covering of $\overline
{B(x,h_1)}$, we can take $j_0$ such that $\langle\xi_{j_0},\xi
\rangle
\geq\gamma/2$. Therefore, $C_{\gamma/2,\xi_{j_0},h_1}(x)\subset
C_{\gamma,\xi,h_1}(x)$. The final statement about the cardinality
follows directly from the construction.\vadjust{\goodbreak}
\end{pf*}

\begin{pf*}{Proof of Theorem~\ref{thunav}}
Note that $F$ is a volume function (and hence a
\textit{Kneser function}) as those considered in Stach{\'o} \cite{sta76}.
According to Lemma~2 in that paper, $F$ is absolutely continuous and
$F^\prime(t)$ exists except for a countable set.
So, there
exist a countable set $N$ and positive constants $s$ and $q$ such that
$F'(t)<q$ $\forall t\in[0,s]\cap N^c$.
If we take $h_2=\min\{s,h_1\}$, where $h_1=\frac{h}{2}\sin(\rho/2)$,
then according with equation~(\ref{equnav}):
\begin{eqnarray*}
\mathbb{E} \bigl(d_{P_X} \bigl(S,\tilde{\mathbb{C}}_{\rho,h}(
\aleph_n) \bigr) \bigr)&\leq& k_2\int_{\{x \in S:d(x,\partial
S)\leq h_2\}}
\sum_{U\in
\mathfrak{U}_{x,\rho,h}} \bigl(1-k_1 \mu(U\cap S)
\bigr)^n \,dx
\\
&&{}+k_2\int_{\{x \in S:d(x,\partial S)> h_2\}} \sum
_{U\in\mathfrak
{U}_{x,\rho,h}} \bigl(1-k_1 \mu(U\cap S)
\bigr)^n \,dx.
\end{eqnarray*}

With respect to the last term, note that $d(x,\partial S)> h_2$ entails
$B(x,h_2)\subset S$ and for all $U\in\mathfrak{U}_{x,\rho,h}$, we
have: $k_1 \mu(U\cap S)\geq k_1 \mu(U\cap B(x,h_2) )= c_0h_2^d$
for some positive $c_0$. Therefore, if $k=k(\rho,h)$ denotes the
cardinality of the set $\mathfrak{U}_{x,\rho,h}$ then
\begin{eqnarray*}
&&\int_{\{x \in S:d(x,\partial S)> h_2\}} \sum_{U\in\mathfrak
{U}_{x,\rho
,h}}
\bigl(1-k_1 \mu(U\cap S) \bigr)^n \,dx
\nonumber
\\[-8pt]
\\[-8pt]
\nonumber
&&\qquad\leq k
\bigl(1-c_0 h_2^d \bigr)^n\mu
\bigl( \bigl\{x \in S\dvtx d(x,\partial S)> h_2 \bigr\} \bigr),
\end{eqnarray*}
which can be upper bounded by $k_3 e^{-nh_2^dc_0}$, for some positive
constant $k_3$.

In order to bound the first integral, note that if
$U=C_{\gamma/2,\xi_j,h_1}(x)$ and $t\leq h_2\leq h_1$ then $U\cap
B(x,t)=C_{\gamma/2,\xi_j,t}(x)$ and so, if $d(x,\partial S)=t$
\[
k_1 \mu(U\cap S)\geq k_1\mu \bigl(U\cap B(x,t)
\bigr)=c_0t^d.
\]
Therefore,
\begin{eqnarray*}
&&\int_{\{x \in S:d(x,\partial S)\leq h_2\}} \sum_{U\in\mathfrak
{U}_{x,\rho,h}}
\bigl(1-k_1 \mu(U\cap S) \bigr)^n \,dx\\
&&\qquad\leq \int
_{\{x
\in S:d(x,\partial S)\leq h_2\}} k \bigl(1-c_0d(x,\partial S)^d
\bigr)^n \,dx
\\
&&\qquad\leq \int_{\{x \in S:d(x,\partial S)\leq h_2\}}
k e^{-c_0nd(x,\partial S)^d} \,dx.
\end{eqnarray*}
Next, let $g(z)= ke^{-c_0nz^d}$. A change of variables leads to
\begin{eqnarray*}
\int_{\{x \in S:d(x,\partial S)\leq h_2\}} k e^{-c_0nd(x,\partial
S)^d} \,dx &=& \int_{\{x \in S:d(x,\partial S)\leq
h_2\}}
g \bigl(d(x,\partial S) \bigr)\,dx\\
&=&\int_{[0,h_2]} g(z)\,dF(z)
\\
&=& \int_{[0,h_2]}g(z)F'(z)\,dz\leq\int
_{[0,h_2]}k_4e^{-c_0nz^d} \,dz,
\end{eqnarray*}
with $k_4$ a positive constant (we have used in the last inequality the
essential boundedness $F'$ in $[0,s]$).
Finally, we have that there exists $k_5>0$ such that
\begin{eqnarray*}
\int_{[0,h_2]}k_4e^{-c_0nz^d} \,dz&=&
n^{-1/d}\int_0^{c_0nh_2^d} k_5e^{-u}u^{(1-d)/d}\,du
\\
&\leq& n^{-1/d}\int_0^{+\infty}
k_5e^{-u}u^{(1-d)/d}\,du\\
& =&\mathcal{O}
\bigl(n^{-1/d} \bigr).
\end{eqnarray*}

Collecting bounds, we get
$\mathbb{E} (d_{P_X} (S,\tilde{\mathbb{C}}_{\rho,h}(\aleph_n)
) )= \mathcal{O} (e^{-nh_2^dc_0 }+n^{-1/d} )=\mathcal
{O}(n^{-1/d})$.
\end{pf*}
\end{appendix}
\section*{Acknowledgements}
We are indebted to Beatriz
Pateiro-L\'opez for insightful corrections. A. Cuevas
is grateful to his colleagues Jes\'us Garc\'{\i}a-Azorero and
Ireneo Peral for useful conversations on the history of Poincar\'e
condition. The constructive comments and criticisms from two
anonymous referees are gratefully acknowledged.

%

%



\printaddresses

\end{document}